\documentclass[leqno]{article}

\usepackage[a4paper,margin=4cm]{geometry}

\usepackage{amsmath,amsrefs,amssymb,amsthm,authblk,bigstrut,comment,enumitem,float,%
leading,longtable,mathrsfs,microtype,
tikz}

\usetikzlibrary{arrows,arrows.meta,calc,decorations,matrix,overlay-beamer-styles,positioning,through}

\tikzset{%
    add/.style args={#1 and #2}{
        to path={%
 ($(\tikztostart)!-#1!(\tikztotarget)$)--($(\tikztotarget)!-#2!(\tikztostart)$)%
  \tikztonodes},add/.default={.2 and .2}}
}

\pgfdeclaredecoration{arrows}{draw}{
\state{draw}[width=\pgfdecoratedinputsegmentlength]{%
  \path [every arrow subpath/.try] \pgfextra{%
    \pgfpathmoveto{\pgfpointdecoratedinputsegmentfirst}%
    \pgfpathlineto{\pgfpointdecoratedinputsegmentlast}%
   };
}}

\tikzset{every arrow subpath/.style={->, draw, thick}}

\newcommand{\foreseven}[9]{
\begin{tikzpicture}[baseline=(a.base),scale=#1,every node/.style={scale=#2},inner sep=0pt,outer sep=0pt]
\node (a) at (1,0) {$#9$};
\node at (2,0) {$#8$};
\node at (3,0) {$#7$};
\node at (4,0) {$#6$};
\node at (4,-1.5) {$#4$};
\node at (5,0) {$#5$};
\node at (6,0) {$#3$};
\end{tikzpicture}
}

\newcommand{\foreight}[9]{
\begin{tikzpicture}[baseline=(a.base),scale=.1,every node/.style={scale=#1},inner sep=0pt,outer sep=0pt]
\node (a) at (1,0) {$#9$};
\node at (2,0) {$#8$};
\node at (3,0) {$#7$};
\node at (4,0) {$#6$};
\node at (5,0) {$#5$};
\node at (5,-1.5) {$#3$};
\node at (6,0) {$#4$};
\node at (7,0) {$#2$};
\end{tikzpicture}
}

\newlength{\halogap}

\newlength\triplesep
\newlength\triplelinewidth
\tikzset{triple/.style={
line width=\triplelinewidth,
black,
 preaction={
  preaction={
    draw,
    line width=2\triplesep+3\triplelinewidth,
    black
   },
   draw,
   line width=2\triplesep+\triplelinewidth,white
  }
 }
}

\tikzset{DynkinNode/.style={circle,minimum size=.7em,inner sep=0pt,font=\scriptsize}}

\newcommand{\GDynkin}[1]{
\begin{tikzpicture}[scale=.8,anchor=base,baseline]
\foreach\kthweight[count=\k] in {#1}{
\ifnum\k=1\node[DynkinNode] (1) at (-1,0) {$\scriptscriptstyle\kthweight$};\fi
\ifnum\k=2\node[DynkinNode] (2) at (-2,0) {$\scriptscriptstyle\kthweight$};\fi
}
\draw[-Implies,double distance=1pt] (2) -- (1);
\draw (2) -- (1);
\end{tikzpicture}
}

\newcommand{\FDynkin}[1]{
\begin{tikzpicture}[scale=.6,anchor=base,baseline]
\foreach\kthweight[count=\k] in {#1}{
\ifnum\k=1\node[DynkinNode] (1) at (-1,0) {$\scriptscriptstyle\kthweight$};\fi
\ifnum\k=2\node[DynkinNode] (4) at (-4,0) {$\scriptscriptstyle\kthweight$};\fi
\ifnum\k=3\node[DynkinNode] (2) at (-2,0) {$\scriptscriptstyle\kthweight$};\fi
\ifnum\k=4\node[DynkinNode] (3) at (-3,0) {$\scriptscriptstyle\kthweight$};\fi
}
\draw[thick] (1) -- (2);
\draw[-Implies,double,thick] (3) -- (2);
\draw[thick] (3) -- (4);
\end{tikzpicture}
}

\newcommand{\EDynkin}[2]{
\begin{tikzpicture}[scale=.4,anchor=base,baseline]
\foreach\kthweight[count=\k] in {#1}{
\pgfmathtruncatemacro{\prevnode}{\k-1}
\ifnum\k=1\node[DynkinNode] (\k) at (0,-1) {$\scriptscriptstyle\kthweight$};\fi 
\ifnum\k>1\node at (4-\k,0) {$\scriptscriptstyle\kthweight$};
\node[DynkinNode] (\k) at (4-\k,0) {$\scriptscriptstyle\kthweight$};
\ifnum\k>2\draw[thick] (\k) -- (\prevnode);\fi
\ifnum\k=4\draw[thick] (\k) -- (1);\fi
\fi
}
\foreach\Node[count=\i] in {#2}
{
\xdef\imax{\i}
\coordinate (AuxNode-\i) at ($(\Node)$);
}
\ifnum\imax=3%
\draw[red,rounded corners] ($(AuxNode-1)+(-\halogap,\halogap)$) -|  ($(AuxNode-2)+(-\halogap,\halogap)$)
    -- ($(AuxNode-2)+(\halogap,\halogap)$) |-
    ($(AuxNode-3)+(\halogap,\halogap)$)
    --($(AuxNode-3)+(\halogap,-\halogap)$) -- ($(AuxNode-1)+(-\halogap,-\halogap)$) -- cycle;
\else\ifnum\imax=2%
\draw[red,rounded corners] ($(AuxNode-1)+(-\halogap,\halogap)$) --  ($(AuxNode-2)+(\halogap,\halogap)$)
    -- ($(AuxNode-2)+(\halogap,-\halogap)$) -- ($(AuxNode-1)+(-\halogap,-\halogap)$)--
    cycle;
\fi
\fi
\end{tikzpicture}
}

\def\rA{{\mathrm A}}
\def\rB{{\mathrm B}}
\def\rC{{\mathrm C}}
\def\rD{{\mathrm D}}
\def\rE{{\mathrm E}}
\def\rF{{\mathrm F}}
\def\rG{{\mathrm G}}

\newcommand\ph\varphi
\renewcommand\le\leqslant
\renewcommand\ge\geqslant
\newcommand\eps\varepsilon
\newcommand\g{{\mathfrak g}}
\newcommand\gl{{\mathfrak{gl}}}
\newcommand\qa{{\mathfrak q}}
\newcommand\la{{\mathfrak l}}

\newcommand\ca{{\mathfrak c}}
\newcommand\z{{\mathfrak z}}
\newcommand\sa{{\mathfrak s}}
\newcommand\sla{{\mathfrak{sl}}}
\newcommand\spa{{\mathfrak{sp}}}
\newcommand\so{{\mathfrak{so}}}
\newcommand\CA{{\mathscr C}}
\newcommand\Ss{{\mathbf S}}
\newcommand\field{{\mathbb F}}
\newcommand{\sslash}{{\mkern-2mu/\mkern-6mu/}}
\newcommand\abr[1]{\left\langle#1\right\rangle}
\newcommand\setof[2]{\left\{#1\mid#2\right\}}
\DeclareMathOperator\ad{ad}
\DeclareMathOperator\rank{rk}

\makeatletter
\let\c@equation=\c@subsection
\def\thmhead#1#2#3{%
  (\thmnumber{
    \@upn{#2}})
   \thmname{#1}%
   \thmnote{ {\the\thm@notefont(#3)}}%
}
\makeatother

\renewcommand{\subsection}[1]{
\refstepcounter{subsection}
\medskip\noindent\textbf{(\thesubsection)\ #1\unskip.}
\ignorespaces}

\renewcommand{\subsubsection}[1]{
\refstepcounter{subsubsection}
\smallskip\noindent\textbf{(\thesubsubsection)}\ \textit{#1\unskip.}
\ignorespaces}

\renewcommand{\thesubsection}{\thesection.\arabic{subsection}}
\renewcommand{\thesubsubsection}{\thesection.\arabic{subsection}.\arabic{subsubsection}}

\numberwithin{equation}{section}

\renewcommand{\theequation}{\thesection.\arabic{equation}}

\newcommand{\theodef}[1]{\newtheorem{#1}[subsection]{#1}}

\theodef{Lemma}
\theodef{Proposition}
\theodef{Theorem}
\theodef{Corollary}
\theodef{Conjecture}

\theoremstyle{definition}

\theodef{Definition}
\theodef{Remark}
\theodef{Example}
\theodef{Examples}

\title
{Semisimple cyclic elements in semisimple Lie algebras}

\author[$^1$]{A. G. Elashvili}
\author[$^1$]{M. Jibladze}
\author[$^2$]{V. G. Kac}

\affil[$^1$]{Razmadze Mathematical Institute, TSU, Tbilisi 0186, Georgia}
\affil[$^2$]{Department of Mathematics, MIT, Cambridge MA 02139, USA}

\date{}

\begin{document}

\maketitle

\begin{abstract}
This paper is a continuation of the theory of cyclic elements in semisimple Lie algebras, developed by Elashvili, Kac and Vinberg. Its main result is the classification of semisimple cyclic elements in semisimple Lie algebras. The importance of this classification stems from the fact that each such element gives rise to an integrable hierarchy of Hamiltonian PDE of Drinfeld-Sokolov type.
\end{abstract}

\section{Introduction}

Let $\g$ be a semisimple finite-dimensional Lie algebra over an algebraically closed field $\field$ of characteristic 0 and let $e$ be a non-zero nilpotent element of $\g$. By the Morozov-Jacobson theorem, the element $e$ can be included in an $\sla(2)$-triple $\sa=\{e,h,f\}$ (unique, up to conjugacy \cite{kostant}), so that $[e,f]=h$, $[h,e]=2e$, $[h,f]=-2f$. Then the eigenspace decomposition of $\g$ with respect to $\ad h$ is a $\mathbb Z$-grading of $\g$:
\begin{equation}\label{grading}
\g=\bigoplus_{j=-d}^d\g_j,\text{ where $\g_{\pm d}\ne0$.}
\end{equation}
The positive integer $d$ is called the \emph{depth} of the nilpotent element $e$.

An element of $\g$ of the form $e+F$, where $F$ is a non-zero element of $\g_{-d}$, is called a \emph{cyclic element}, associated to $e$. In \cite{kostant} Kostant proved that any cyclic element, associated to a principal ($=$ regular) nilpotent element $e$, is regular semisimple, and in \cite{springer} Springer proved that any cyclic element, associated to a subregular nilpotent element of a simple exceptional Lie algebra, is regular semisimple as well, and, moreover, found two more distinguished nilpotent elements in E$_8$ with the same property.

A non-zero nilpotent element $e$ of $\g$ is called of \emph{nilpotent} (resp. \emph{semisimple}) \emph{type} if all cyclic elements, associated to $e$, are nilpotent (resp. there exists a semisimple cyclic element, associated to $e$). If neither of the above cases occurs, the element $e$ is called of \emph{mixed type} \cite{ekv}.

It is explained in the introduction to \cite{ekv} how to reduce the study of cyclic elements to the case when $\g$ is simple. Therefore, we shall assume from now on that $\g$ is simple, unless otherwise stated.

An important r\^ole in the study of cyclic elements, associated to a non-zero nilpotent element $e$, is played by the centralizer $\z(\sa)$ in $\g$ of the $\sla(2)$-triple $\sa$ and by its centralizer $Z(\sa)$ in the connected adjoint group $G$. Since $h\in\sa$, the group $Z(\sa)$ preserves the grading \eqref{grading}.

Let us state now some of the main results from \cite{ekv}.

\begin{Theorem}\label{odepth}
A nilpotent element $e$ is of nilpotent type iff the depth $d$ of $e$ is odd. In this case the group $Z(\sa)$ has finitely many orbits in $\g_{-d}$, hence zero is the only closed orbit.
\end{Theorem}

\begin{Theorem}\label{edepth}
If a non-zero nilpotent element $e$ has even depth, then the representation of $Z(\sa)$ in $\g_{-d}$ is orthogonal, i.~e. preserves a non-degenerate invariant symmetric bilinear form $(\cdot,\cdot)$. Consequently, by \cite{L} the union of closed orbits of $Z(\sa)$ in $\g_{-d}$ contains a non-empty Zariski open subset.
\end{Theorem}

Let
\begin{equation}
\Ss_\g(e)=\setof{F\in\g_{-d}}{\text{$e+F$ is semisimple in $\g$}}.
\end{equation}

\begin{Theorem}\label{aboutsg}
Let $e\in\g$ be a nilpotent element of semisimple type. Then
\begin{itemize}
\item[(a)] $\Ss_\g(e)$ contains a non-empty Zariski open subset in $\g_{-d}$.
\item[(b)] If $F\in\Ss_\g(e)$, then the $Z(\sa)$-orbit of $F$ in $\g_{-d}$ is closed.
\end{itemize}
\end{Theorem}

Thus, $\Ss_\g(e)$ consists of closed $Z(\sa)$-orbits in $\g_{-d}$, and in order to classify semisimple cyclic elements, we need to describe, for each nilpotent element $e$ of semisimple type, the complement to $\Ss_\g(e)$ in $\g_{-d}$, which we call the \emph{singular subset of $\g_{-d}$}.

Recall that the dimension of $\g_{-d}\sslash Z(\sa)$ is called the \emph{rank} of the nilpotent element $e$, and is denoted by $\rank e$.

The representation of the group $Z^\circ(\sa)$, the unity component of $Z(\sa)$, in $\g_{-d}$ is given in \cite{ekv} for each nilpotent element $e$, whose type is not nilpotent. It follows from this description that all these representations are strongly polar in the following sense (see Section \ref{polarsection} for details). We call a representation of a reductive group $S$ in a vector space $V$ \emph{strongly polar} if it is polar in the sense of \cite{dadokac}, and every maximal subspace of $V$, consisting of vectors with closed $S$-orbits, called a \emph{Cartan subspace}, has dimension equal to that of $V\sslash S$ and all Cartan subspaces in $V$ are conjugate by $S$. (Recall that $V\sslash S:=\operatorname{Specm}\field[V]^S$.) This is a stronger version of the definition of a polar representation, introduced in \cite{dadokac}, but it is conjectured there that these definitions are equivalent.

Note that, by definition, $\rank e$ is equal to the dimension of a Cartan subspace for $Z(\sa)$ in $\g_{-d}$.

The basic notion of the theory of cyclic elements is that of a \emph{reducing} subalgebra, which we give here for nilpotent elements of semisimple type.

\begin{Definition}
Let $e$ be a nilpotent element of semisimple type in $\g$. A subalgebra $\qa$ of $\g$ is called a \emph{reducing subalgebra} for $e$ if $\qa$ is semisimple, contains $\sa$, hence $\ad h$ induces $\mathbb Z$-grading $\qa=\bigoplus_j\qa_j$, and $Z(\sa)(\qa_{-d})$ contains a Zariski open subset in $\g_{-d}$.
\end{Definition}

The first result of the paper, presented in Section \ref{sectionred}, is the following theorem, which is a stronger version for elements of semisimple type of Theorem 3.14 from \cite{ekv}.

\begin{Theorem}\label{redcartan}
If $e$ is a nilpotent element of semisimple type in $\g$, then there exists a reducing subalgebra $\qa$ for $e$, such that $\qa_{-d}$ is a Cartan subspace of the representation of $Z(\sa)$ in $\g_{-d}$.
\end{Theorem}

Unfortunately, we do not know a proof of this theorem without a case-wise verification using Tables \ref{reductions}, \ref{reductionsf4}, \ref{reductionse6}, \ref{reductionse7}, \ref{reductionse8} and \ref{irreds}. It turns out that the minimal Levi subalgebra, containing $e$, does the job for most of the cases. This fails only for one kind of nilpotent elements in $\g$ for each of the types B$_n$, C$_n$ and F$_4$.

Using \eqref{redzs} below, Theorem \ref{redcartan} reduces the classification of semisimple cyclic elements, associated to a non-zero nilpotent element $e$, to the case when $e$ is a distinguished nilpotent element in $\g$, namely when the group $Z(\sa)$ is finite. Obviously we may assume in addition that $\g$ does not contain a smaller reducing subalgebra for $e$. In this case the nilpotent element $e$ of semisimple type is called \emph{irreducible}.

Note that, obviously, $\Ss_\g(e)$ is $Z(\sa)$-invariant, hence conical (see Proposition \ref{propcon} below). In particular, if $e$ is a semisimple type nilpotent element of rank 1, taking $F_0\in\g_{-d}$, such that $(F_0,F_0)\ne0$, we obtain (using Theorem \ref{edepth})
$$
\Ss_\g(e)=\field^*Z(\sa)F_0=\setof{F\in\g_{-d}}{(F,F)\ne0}.
$$
It turns out that there are very few irreducible nilpotent elements of rank $>1$ in simple Lie algebras: one of rank 2 in $\so(4k)$ and F$_4$, one of rank 3 in E$_7$, two of rank 2 in E$_8$ and one of rank 4 in E$_8$. These cases are treated in Section \ref{irredsection}, giving thereby a complete description of the set $\Ss_\g(e)$ for all simple Lie algebras $\g$ and nilpotent elements $e$ of semisimple type.

Namely, an arbitrary nilpotent element $e$ of semisimple type in a simple Lie algebra $\g$ is irreducible in a direct sum of simple Lie algebras $\qa_1\oplus\cdots\oplus\qa_s$, containing $e$ with non-zero projections $e_i$ to $\qa_i$ of the same depth as $e$ in $\g$, such that $e_i$ is irreducible in $\qa_i$ and $e+F$ with $F=F_1+...+F_s\in\qa_1\oplus\cdots\oplus\qa_s$ is semisimple iff each $e_i+F_i$ is semisimple.

In the last Section \ref{algsection} we relate the problem of finding all semisimple cyclic elements, associated to a nilpotent element $e$ of depth $d$, to an algebra structure on the subspace $\g_{-d}$, defined by the formula (recall that $d$ is even if $e$ is not of nilpotent type \cite{ekv})
\begin{equation}\label{mdprod}
x*y=[(\ad e)^{\frac d2}x,y],\quad x,y\in\g_{-d}.
\end{equation}
One easily shows that this product is commutative (resp. anticommutative) if $\frac12d$ is odd (resp. even).

It is well known that an even nilpotent element $e$ of depth $d=2$ is always of semisimple type, and the product \eqref{mdprod} defines on $\g_{-2}$ a structure of a simple Jordan algebra (in fact, all simple Jordan algebras are thus obtained \cite{J}).

It turns out that for an irreducible nilpotent element $e$ of rank ($=\dim\g_{-d}:=n$) $\ge2$ the algebra \eqref{mdprod} is always a commutative algebra, denoted by $\CA_\lambda(n)$, for some particular $\lambda\in\field$, which in a basis $p_1$, ..., $p_n$ has multiplication table
\begin{equation}\label{multab}
p_i^2=p_i,\quad\text{$p_ip_j=\lambda(p_i+p_j)$ if $i\ne j$.}
\end{equation}
For $\lambda\ne\frac12$ the algebra $\CA_\lambda(n)$ has $\le2^n-1$ nonzero idempotents, in fact, except for an easily describable finite set of exceptions, exactly $2^n-1$ of them. In particular this is so in all cases that occur in our situations. For example, if $n=2$, then the values of $\lambda$ are as follows:

$\g=\so(4k)$: $\lambda=1-k$; $\g=$F$_4$: $\lambda=-\frac13$; $\g=$E$_8$: $\lambda=-1$ and $-\frac27$.

We compute the algebra \eqref{mdprod} for all nilpotent elements of semisimple type. Obviously this algebra is the same as for the corresponding irreducible nilpotent element in the cases when $e$ is such that $\rank e=\dim\g_{-d}$. Remarkably, it turns out that in all other cases this algebra is either a direct sum of at most two simple Malcev algebras (including the 1-dimensional Lie algebra), which happens iff $\frac d2$ is even, or a simple Jordan algebra, which happens iff $\frac d2$ is odd.

What does it have to do with the main problem in question? It turns out that one can describe the singular subset $\g_{-d}\setminus\Ss_\g(e)$ in terms of this algebra. We show that for an irreducible nilpotent element $e$ of depth $d$ with odd $\frac d2$ the singular subset consists of those $F\in\g_{-d}$, which are contained in a proper subalgebra of the algebra \eqref{mdprod}. For example, in the case $n=2$ the singular subset consists of scalar multiples of the three non-zero idempotents (see \eqref{multab}):
$$
p_1,p_2,\text{ and $\frac{p_1+p_2}{2\lambda+1}$}.
$$
In general, for $n=3$ and $n=4$, the singular subset consists of the union of spans of $n-1$ linearly independent idempotents, namely, it is a union of $\frac{n(n+1)}2$ hyperplanes in the $n$-dimensional space.

For an arbitrary nilpotent element $e$ of semisimple type either there is a reducing subalgebra which is a direct sum of isomorphic simple Lie algebras with each projection of $e$ to them being a nilpotent element of rank 1 (in fact, principal), or the depth $d$ is such that $\frac d2$ is odd. In the latter case the algebra $\g_{-d}$ with product \eqref{mdprod} is commutative and its Cartan subspace is a subalgebra $\ca$, isomorphic to one of the algebras corresponding to irreducible nilpotent elements. Then the singular subset for $e$ is equal to $Z(\sa)(\ca_{\mathrm{sing}})$, where $\ca_{\mathrm{sing}}$ is the singular subset of $\ca$ (described above).

We list in Table \ref{irreds} (see Section \ref{irredsection}) all irreducible nilpotent elements of semisimple type in all simple Lie algebras, and in Tables \ref{reductions}, \ref{reductionse6}, \ref{reductionse7}, \ref{reductionse8}, \ref{reductionsf4} (see Section \ref{sectionred}) all non-irreducible nilpotent elements of semisimple type in simple Lie algebras of types A, B, C, D; E$_6$; E$_7$; E$_8$; F$_4$ and G$_2$ respectively (using the tables in \cite{ekv}), along with their depth, rank, the minimal reducing subalgebra $\qa_{\min}$ (by its number in Table \ref{irreds}), and the structure of the algebra $(\g_{-d},*)$.

Many of our results are proved in the tradition of ancient Greeks: look at the tables! It would be interesting to find unified proofs of such claims. Here are some of them:
\begin{itemize}
\item[(a)] If $d$ is odd, then the linear group $Z(\sa)|\g_{-d}$ is Sp$(n)$ (we know a priori that this is a subgroup of Sp$(n)$ with finitely many orbits).
\item[(b)] If $d$ is even, then the linear group $Z(\sa)|\g_{-d}$ is strongly polar and $\g_{-d}$ is a sum of at most two irreducible modules.
\item[(c)] If $d$ is divisible by $4$, then $(\g_{-d},*)$ is a Malcev algebra.
\item[(d)] If $s:=\dim\g_{-d}>1$ and the group $Z(\sa)|\g_{-d}$ is finite, then the algebra $(\g_{-d},*)$ has exactly $2^s$ idempotents and the singular set is a union of hyperplanes, spanned by idempotents, their number being $s(s+1)/2$.
\item[(e)] If $e$ is of semisimple type, the group $Z(\sa)|\g_{-d}$ is infinite, and $d/2$ is odd, then $(\g_{-d},*)$ is a simple Jordan algebra.
\end{itemize}

In conclusion of the introduction recall that one of the applications of the study of semisimple cyclic elements is that to regular elements in Weyl groups \cites{kostant,springer,ekv}. Another application goes back to the work of Drinfeld and Sokolov \cite{drisok}, where they used the principal cyclic elements of simple Lie algebras to construct integrable Hamiltonian hierarchies of PDE of KdV type (the KdV arising from $\sla(2)$). This work followed by series of papers by various authors, where the method of \cite{drisok} was extended to other semisimple cyclic elements. In complete generality this has been done in \cite{dSkv}, where to each semisimple cyclic element, considered up to a non-zero constant factor and up to conjugacy by $Z(\sa)$, an integrable Hamiltonian hierarchy of PDE was constructed.

The contents of the paper is as follows. After explaining the basic notions, the goal, and the motivations of the paper in the Introduction, we discuss the notions of polar and strongly polar linear reductive algebraic groups in Section 2 (Theorems \ref{cartanweyl} and \ref{thepolar}). The reason for it is Proposition \ref{proporpolar}, which claims that the linear group $Z(\sa)|\g_{-d}$ is strongly polar. This, along with Theorem \ref{aboutsg}, restricts considerably the possibilities for $F\in\g_{-d}$, such that the cyclic element $e+F$ is semisimple.

In Section \ref{irredsection} we list irreducible nilpotent elements $e$ of semisimple type in Table \ref{irreds}. By definition, they don't admit a nontrivial reducing subalgebra, and consequently the group $Z(\sa)|\g_{-d}$ is finite (these finite linear groups are listed in Table \ref{irreds}). Theorem \ref{conj} describes an explicit parametrization of the set $\Ss_\g(e)$ for all nilpotent elements $e$ from Table \ref{irreds} in simple Lie algebras $\g$.

In Section \ref{sectionred} for each nilpotent element $e$ of semisimple type in a simple Lie algebra $\g$ we exhibit a (semisimple) reducing subalgebra where $e$ is irreducible. This reduces the description of the set $\Ss_\g(e)$ to the irreducible nilpotent elements of semisimple type from Table \ref{irreds}. The obtained information on nilpotents $e$ of semisimple type in simple classical Lie algebras, in $\g$ of type F and G, and in $\g$ of type E$_6$, of type E$_7$, and of type E$_8$, is given in Tables \ref{reductions}, \ref{reductionsf4}, \ref{reductionse6}, \ref{reductionse7}, and \ref{reductionse8}, respectively.

Finally, in Section \ref{algsection} we study the algebra $(\g_{-d},*)$, associated to a nilpotent $e$ of semisimple type by formula \eqref{mdprod}. It is a generalization of the well-known construction of simple Jordan algebras when $d=2$. These algebras are explicitly described by Theorem \ref{algm}. In Theorem \ref{theoS} we provide the description of the set $\Ss_\g(e)$ in terms of these algebras.

We added to the paper three appendices. In Appendix A we describe for each odd nilpotent element $e$ the even subalgebra $\g^{\mathrm{ev}}=\bigoplus_{j\in\mathbb Z}\g_{2j}$. In Appendix B we describe the algebras $(\g_{-d},*)$ for all nilpotent elements of mixed type in $\g$. In Appendix C we describe chains for all nilpotent elements in $\g$, which is a generalization of the decomposition into unions of Jordan blocks of the same size in $\g=\sla(n)$.

Throughout the paper the base field $\field$ is an algebraically closed field of characteristic zero.

We are grateful to E. B. Vinberg for numerous discussions and suggestions, and to a referee for a large number of questions and corrections. All the calculations were made possible thanks to the \texttt{GAP} system for computational algebra, and especially the \texttt{GAP} package \texttt{SLA} by Willem de Graaf \cite{deGraaf}, who also provided several helpful emails explaining its usage. The paper was completed while all three authors visited, in the summer of 2019, the IHES, France, whose hospitality is gratefully acknowledeged.

\section{Polar representations and reducing subalgebras}\label{polarsection}

Let $G$ be a reductive subgroup of GL$(V)$, where $V$ is a finite-dimensional vector space over $\field$, and let $\g\subseteq\gl(V)$ be its Lie algebra. Let $v\in V$ be such that its orbit $G(v)$ is closed. Let
\begin{equation}\label{cv}
\ca_v=\setof{x\in V}{\g(x)\subseteq\g(v)}.
\end{equation}
Then \cite{dadokac} $\dim\ca_v\le\dim V\sslash G$. The linear reductive group $G$ is called \emph{polar} if
\begin{equation}\label{defpolar}
\dim\ca_v=\dim V\sslash G,
\end{equation}
and in this case $\ca_v$ is called a \emph{Cartan subspace} of $V$. Note that, by definition, $G$ is polar iff its identity component is.

The following results are either proved in \cite{dadokac} or easily follow from it.

Let $\ca\subset V$ be a Cartan subspace, and let $N(\ca)=\setof{g\in G}{g(\ca)=\ca}$, $Z(\ca)=\setof{g\in G}{\text{$g(v)=v$ for all $v\in\ca$}}$. Then $N(\ca)/Z(\ca)$ is called the \emph{Weyl group} of the polar linear group $G$.

\ 

\ 

\begin{Theorem}\label{cartanweyl}
Let $G\subset\operatorname{GL}(V)$ be a polar linear group, let $\ca\subset V$ be a Cartan subspace, and let $W\subset\operatorname{GL}(\ca)$ be the Weyl group of $\ca$. Then
\begin{itemize}
\item[(a)] Any Cartan subspace $\ca_1\subset V$ is conjugate by $G$ to $\ca$.
\item[(b)] The Weyl group $W$ is finite and any closed orbit of $G$ intersects $\ca$ by an orbit of $W$. Furthermore, ${\mathbb C}[V]^G\xrightarrow\sim{\mathbb C}[\ca]^W$ via restriction.
\item[(c)] If $G$ is connected, then the Weyl group $W$ is generated by unitary reflections. If $G$ is orthogonal, then $G\cdot\ca$ is Zariski dense in $V$ and $W$ is generated by orthogonal reflections.
\end{itemize}
\end{Theorem}

\begin{proof}
Claim (a) is a part of Theorem 2.3 from \cite{dadokac}.

Claim (b) is Lemma 2.7 and Theorems 2.8, 2.9 from \cite{dadokac}.

Claim (c), except for the second part, is Theorem 2.10 from \cite{dadokac}.

If $V$ is orthogonal, i.~e. has a non-degenerate symmetric $G$-invariant bilinear form $(\cdot,\cdot)$, then the generic $G$-orbit is closed by \cite{L}, hence the restriction of $(\cdot,\cdot)$ to $\ca$ is non-degenerate $W$-invariant, hence the reflections in $W$ are orthogonal.
\end{proof}

\

\begin{Theorem}\label{thepolar}\
\begin{itemize}
\item[(a)] A direct sum of linear reductive groups $G_i\subset\operatorname{GL}(V_i)$ is polar iff all $G_i\subset\operatorname{GL}(V_i)$ are polar.
\item[(b)] If $G\subseteq\operatorname{GL}(V)$ is a reductive subgroup and $\dim V\sslash G\le1$ or $=\dim V_0$, where $V_0$ is the zero weight space for $\g$ in $V$, then $G$ is polar, $V_0$ being a Cartan subspace in the second case.
\item[(c)] All theta-groups are polar.
\item[(d)] For a theta-group, any subspace $\ca\subset V$ consisting of semisimple elements, and such that $\dim\ca=\dim V\sslash G$, is a Cartan subspace. Consequently all theta-groups are strongly polar.
\end{itemize}
\end{Theorem}

\begin{proof}
Claims (a) and (b) are obvious.

Claim (c) was stated without proof in \cite{dadokac}. It follows easily from \cite{V}. Indeed, recall \cites{k75,V} that a theta-group is obtained by considering the grading defined by an order $m$ automorphism $\theta$ of a reductive Lie algebra $\mathfrak p$:
\begin{equation}\label{decomp}
\mathfrak p=\bigoplus_{j\in\mathbb Z/m\mathbb Z}\mathfrak p_j.
\end{equation}
Then the connected linear algebraic group $P_0$ with Lie algebra $\mathfrak p_0$, acting on $\mathfrak p_1$, is called a theta group. It was proved in \cite{V} that if $\ca\subset\mathfrak p_1$ is a maximal abelian subalgebra, consisting of semisimple elements, then
\begin{equation}\label{dimca}
\dim\ca=\dim\mathfrak p_1\sslash G_0.
\end{equation}
Consider the weight space decomposition of $\mathfrak p$ with respect to $\ca$: $\mathfrak p=\bigoplus_{\lambda\in\ca^*}\mathfrak p_\lambda$, so that $\mathfrak p_{\lambda=0}$ is the centralizer of $\ca$ in $\mathfrak p$. Take $v\in\ca$, such that $\lambda(v)\ne0$ for all $\lambda\ne0$ such that $\mathfrak p_\lambda\ne0$. Then, obviously, $[\mathfrak p,x]\subseteq[\mathfrak p,v]$ for $x\in\ca$. Considering the projection of $\mathfrak p$ to $\mathfrak p_1$ with respect to \eqref{decomp}, we deduce that $[\mathfrak p_0,x]\subseteq[\mathfrak p_0,v]$, which together with \eqref{dimca} shows that $\ca$ is a Cartan subspace, proving (c).

Finally claim (d) follows from \cite{MTT}, as claimed in \cite{dadokac}. Indeed if $\mathfrak p_1\subset\mathfrak p$ is as in \eqref{decomp} and if $\ca\subset\mathfrak p_1\subset\mathfrak p\subset\operatorname{End}\mathfrak p$ is a subspace, consisting of semisimple elements, then, by \cite{MTT} it is abelian. Hence, if, in addition, \eqref{dimca} holds, $\ca$ is a maximal abelian subalgebra in $\mathfrak p_1$, consisting of semisimple elements. Therefore, by the discussion proving (c), it is a Cartan subspace.
\end{proof}

\begin{Remark}
As D. Panyushev pointed out to the third author, the group SL(2), acting on the direct sum $V$ of the 2- and 3-dimensional irreducible representations, is not polar, though it has a 2-dimensional subspace consisting of elements with a closed orbit and $\dim V\sslash$SL$(2)=2$.
\end{Remark}

Examples of orthogonal theta-groups:
\begin{itemize}
\item[1)] adjoint representations,
\item[2)] nontrivial representations of F$_4$ and G$_2$ of minimal dimension,
\item[3)] standard representation of SO$_n$,
\item[4)] symmetric square of the standard representation of SO$_n$,
\item[5)] skew-symmetric square of the standard representation of Sp$_n$.
\end{itemize}

\begin{Proposition}\label{proporpolar}
If $e$ is a nilpotent element of semisimple or mixed type in a simple Lie algebra $\g$, then the image of the representation of $Z(\sa)$ in $\g_{-d}$ is orthogonal polar. Moreover any of its subspaces $\ca$ of dimension equal to $\dim\g_{-d}\sslash Z(\sa)$ consisting of elements with closed orbits is a Cartan subspace. Consequently the linear reductive group $Z(\sa)|\g_{-d}$ is strongly polar.
\end{Proposition}

\begin{proof}
The first claim is Theorem \ref{edepth} (by Theorem \ref{odepth}). Just a look at Tables 5.1--5.4 from \cite{ekv} (cf. Tables \ref{reductions}, \ref{reductionsf4}, \ref{reductionse6}, \ref{reductionse7}, \ref{reductionse8} below for semisimple type nilpotent elements) shows that the linear reductive group in question is a direct sum of theta-groups (and 1-dimensional trivial linear groups), see Examples. Hence the proposition follows from Theorem \ref{thepolar} (d).
\end{proof}

\begin{Remark}
It follows from \cite{ekv}, Lemma 1.2, that if $e$ is of nilpotent type, then $\dim\g_{-d}\sslash Z(\sa)=0$, consequently the image of the representation of $Z(\sa)$ in $\g_{-d}$ is polar as well. In fact nilpotent elements of nilpotent type exist in case of classical simple Lie algebras only in $\g=\so(N)$, $N\ge7$, and those correspond to the partition $[2m+1>\underbrace{2m=\cdots=2m}_{\text{$2k$ times}}>\cdots]$ \cite{ekv}, in which case the image of the representation of $Z(\sa)$ in $\g_{-d}$ is the standard representation of Sp$(2k)$. Nilpotent elements of nilpotent type in exceptional Lie algebras are listed in \cite{ekv}*{Table 1.1}. One can show that for all of them the image of the representation of $Z(\sa)$ in $\g_{-d}$ is again the standard representation of Sp$(2k)$ for some $k$. Furthermore, this $k$ equals 1 in all cases, with the following four exceptions:
\begin{itemize}
\item[] $\g=$E$_7$, $e=4$A$_1$, $k=3$;
\item[] $\g=$E$_8$, $e=4$A$_1$, $k=4$; $e=2($A$_2+$A$_1)$, $k=2$; $e=2$A$_3$, $k=2$.
\end{itemize}
\end{Remark}

\begin{Definition}
A semisimple subalgebra $\mathfrak q$ of $\g$ is called \emph{reducing} for a nilpotent element $e$ of semisimple type, if $\mathfrak q$ contains $\sa$ (hence is a $\mathbb Z$-graded subalgebra) and a Cartan subspace of $Z_{\mathfrak q}(\sa)$ in $\mathfrak q_{-d}$ is a Cartan subspace of $Z_\g(\sa)$ in $\g_{-d}$.
\end{Definition}

It follows from Proposition \ref{proporpolar} that in the case of $e$ of semisimple type this definition is equivalent to that in \cite{ekv}. Moreover, the following is an easy consequence of results in \cite{ekv}*{Section 3}:

\begin{Proposition}
If a nilpotent element $e\in\g$ is of semisimple type, then a semisimple subalgebra $\qa$ of $\g$ is reducing for $e$ if and only if it contains $e$ and $e$ has the same depth and rank in $\qa$ as in $\g$.
\end{Proposition}

\begin{Example}
Let $e$ be a nilpotent element of semisimple type in $\g$. Let $\qa_{\max}$ denote the subalgebra of $\g$, generated by $e$ and $\g_{-d}$. It follows from \cite{ekv}*{Theorem 3.3 and Propositions 3.9, 3.10} that $\qa_{\max}$ is a reducing subalgebra for $e$ in $\g$. Note that the derived subalgebra of $\g^{\mathrm{ev}}=\bigoplus_{j\in\mathbb Z}\g_{2j}$ is a reducing subalgebra for $e$, which might be larger than $\qa_{\max}$, so this notation is misleading. One may think of $\qa_{\max}$ as the maximal useful reducing subalgebra.
\end{Example}

\begin{Example}
Let $\sigma$ be a diagram automorphism of $\g$ and $e\in\g$ a $\sigma$-invariant nilpotent element of semisimple type, such that $\g_{-d}^\sigma=\g_{-d}$. Then $\g^\sigma$ is a reducing subalgebra for $e$. This happens if $e$ is a principal nilpotent element in $\g=\sla(2k)$, $\so(2k)$ or E$_6$ and order$(\sigma)=2$, or in $\g=\so(8)$ and order$(\sigma)=3$.
\end{Example}

Recall that the \emph{rank} $\rank_\g(e)$ of $e$ in $\g$ is the dimension of $\g_{-d}\sslash Z(\sa)$. Note that for any reducing subalgebra $\qa$ of a nilpotent element $e$ of semisimple type in $\g$ we have, in view of Theorem \ref{redcartan},
\begin{equation}\label{redzs}
\Ss_\g(e)=Z(\sa)\Ss_\qa(e).
\end{equation}
Indeed, according to Theorem \ref{aboutsg} (b) if $e+F$ is a semisimple element of $\g$ (resp. $\qa$), then the $Z(\sa)$-orbit of $F$ in $\g_{-d}$ (resp. $\qa_{-d}$) is closed.
So, since both representations are polar, we may assume that $F$ lies in a Cartan subspace of $\qa_{-d}$, which is a Cartan subspace of $\g_{-d}$, since $\qa$ is a reducing subalgebra. Thus we can reduce description of $\Ss_\g(e)$ to that for $\Ss_\qa(e)$.

\begin{Proposition}\label{propcon}
The set $\Ss_\g(e)\subset\g_{-d}$ is conical, i.~e. if $F\in\Ss_\g(e)$, then $cF\in\Ss_\g(e)$ for any $c\in\field\setminus\{0\}$.
\end{Proposition}
\begin{proof}
Let $\lambda(t)\subset G$ be the 1-parameter subgroup, corresponding to $h\in\g$ from $\sa$. Then
$$
\lambda(t)(e+F)=t^2e+t^{-d}F,\quad t\in\field\setminus\{0\},
$$
hence $e+t^{-d-2}F$ lies in $\Ss_\g(e)$ if $F$ does.
\end{proof}

\section{Irreducible nilpotent elements of semisimple type}\label{irredsection}

Recall that a nilpotent element $e$ of semisimple type in a simple Lie algebra $\g$ is called irreducible if it does not admit a nontrivial reducing subalgebra different from $\g$ \cite{ekv}. Irreducible nilpotent elements are listed in Table \ref{irreds} below (where $k\ge1$). Recall that in all these cases the linear group $Z(\sa)|\g_{-d}$ is finite and $\dim\g_{-d}=\rank e$. It turns out, using \cite{ale}, \cite{CM}*{Corollary 6.1.6}, that in all cases this finite group is S$_n$ for $n=1,2,3,5$.

Action of this group, as well as the actions of the component groups of $Z(\sa)$ on $\g_{-d}$ in general, are computed in the following way.
First, using the SLA command $\mathtt{FiniteOrderInnerAutomorphisms}$ \cite{deGraaf}, one finds those inner automorphisms of $\g$ of required orders which fix a minimal regular semisimple subalgebra containing $e$. That command provides Kac diagrams of these automorphisms \cite{OV}*{p. 213}; from the Kac diagrams one determines actions of these automorphisms on $\g_{-d}$. In this way we find that if $Z(\sa)|\g_{-d}=$ S$_n$ and $\dim\g_{-d}=n$ (resp. $=n-1$), the group S$_n$ acts on $\g_{-d}$ as the permutation representation (resp. the nontrivial $n-1$-dimensional irreducible representation). We denote the latter by $\sigma_n$, so that the former is $\sigma_n\oplus{\bf1}$. In the last column we list the structure of the algebra $(\g_{-d},*)$; the symbol $1$ there stands for a $1$-dimensional algebra with non-zero (resp.~zero) multiplication if $d/2$ is odd (resp.~even). The algebras $\CA_\lambda(n)$ are defined by \eqref{multab}.

The irreducible nilpotents of semisimple type are listed in the following table (where $k\ge1$):

\begin{table}[H]\renewcommand\thetable{1}\caption{Irreducible nilpotent elements of semisimple type\label{irreds}}

\

\hskip2em\resizebox{.9\textwidth}{!}{%
\begin{tabular}{r|l|r|r|r|r|r}
\#&$\g$&nilpotent $e$&depth&rank&$Z(\mathfrak s)|\g_{-d}$&$(\g_{-d},*)$\\
\hline\hline
1$_k$&$\sla(2k+1)$&$[2k+1]$&$4k$&1\bigstrut[t]&{\bf1}&{\bf1}\\
\hline
2$_k$&$\spa(2k)$&$[2k]$&$4k-2$&1\bigstrut[t]&{\bf1}&{\bf1}\\
\hline
3$_k$&$\so(2k+1)$, $k\ne3$&$[2k+1]$&$4k-2$&1\bigstrut[t]&{\bf1}&{\bf1}\\
4$_k$&$\so(4k+4)$&$[2k+3,2k+1]$&$4k+2$&2&{\bf1} $\oplus$ {\bf1}&$\CA_{-k}(2)$\\
\hline
5&G$_2$&G$_2$\hfill\GDynkin{2,2}&10&1&{\bf1}&{\bf1}\bigstrut[t]\\
\hline
6&F$_4$&F$_4\hfill\FDynkin{2,2,2,2}$&22&1&{\bf1}&{\bf1}\bigstrut[t]\\
7&F$_4$&F$_4(a_2)$\hfill\FDynkin{2,0,0,2}&10&2&$\sigma_2\oplus{\bf1}$&$\CA_{-\frac13}(2)$\bigstrut[b]\\
\hline
8&E$_6$&E$_6(a_1)$\hfill\EDynkin{2,2,2,0,2,2}{1}&16&1&{\bf1}&{\bf1}\bigstrut[t]\\
\hline
9&E$_7$&E$_7$\hfill\EDynkin{2,2,2,2,2,2,2}{1}&34&1&{\bf1}&{\bf1}\bigstrut[t]\\
10&E$_7$&E$_7(a_1)$\hfill\EDynkin{2,2,2,0,2,2,2}{1}&26&1&{\bf1}&{\bf1}\\
11&E$_7$&E$_7(a_5)$\hfill\EDynkin{0,0,0,2,0,0,2}{1}&10&3&$\sigma_3\oplus{\bf1}$&$\CA_{-\frac13}(3)$\bigstrut[b]\\
\hline
12&E$_8$&E$_8$\hfill\EDynkin{2,2,2,2,2,2,2,2}{1}&58&1&{\bf1}&{\bf1}\bigstrut[t]\\
13&E$_8$&E$_8(a_1)$\hfill\EDynkin{2,2,2,0,2,2,2,2}{1}&46&1&{\bf1}&{\bf1}\\
14&E$_8$&E$_8(a_2)$\hfill\EDynkin{2,2,2,0,2,0,2,2}{1}&38&1&{\bf1}&{\bf1}\\
15&E$_8$&E$_8(a_4)$\hfill\EDynkin{0,2,0,2,0,2,0,2}{1}&28&1&{\bf1}&{\bf1}\\
16&E$_8$&E$_8(a_5)$\hfill\EDynkin{0,2,0,2,0,0,2,0}{1}&22&2&$\sigma_2\oplus{\bf1}$&$\CA_{-\frac27}(2)$\\
17&E$_8$&E$_8(a_6)$\hfill\EDynkin{0,0,0,2,0,0,2,0}{1}&18&2&$\sigma_3$&$\CA_{-1}(2)$\\
18&E$_8$&E$_8(a_7)$\hfill\EDynkin{0,0,0,0,2,0,0,0}{1}&10&4&$\sigma_5$&$\CA_{-\frac13}(4)$
\end{tabular}%
}
\end{table}

Irreducibility follows from the fact that in any reducing subalgebra the nilpotent $e$ must have the same depth and rank. For the rank 1 case, dimension of $\g_{-d}$ is 1, and $e$ together with any nonzero $F\in\g_{-d}$ generates $\g$ as an algebra. Now any reducing subalgebra must contain $e$ and a scalar multiple of $F$, so must coincide with $\g$. For rank 2, examining all pairs of cases with equal depth it turns out that none of them can be embedded into each other. There is only one case of rank 3 and only one of rank 4, which implies irreducibility for these ranks.

In this, as well as in all subsequent tables, a nilpotent element is represented by the corresponding partition for classical types, and by its label and the weighted Dynkin diagram for exceptional types. For the latter we use labels from \cite{CM}.

We will describe explicitly in the next Section the minimal reducing subalgebra $\qa_{\min}$ for each nilpotent element $e$ of semisimple type, where it is irreducible. We list there in Tables \ref{reductions}, \ref{reductionse6}, \ref{reductionse7}, \ref{reductionse8}, \ref{reductionsf4} all reducible nilpotent elements of semisimple type and their minimal reducing subalgebras $\qa_{\min}$ by their number $1_k$ -- $4_k$, $5$ -- $18$ from Table \ref{irreds}.

We now turn to the description of the sets $\Ss_\g(e)\subset\g_{-d}$ for irreducible nilpotent elements. The following theorem has been checked with the aid of computer.

\begin{Theorem}\label{conj}
For every irreducible nilpotent element $e$ of semisimple type in a simple Lie algebra $\g$ with $\dim\g_{-d}=m$ there exists an explicit linear isomorphism
$$
I:\g_{-d}\cong\setof{(z_0,...,z_m)\in\field^{m+1}}{z_0+...+z_m=0}
$$
such that
$$
(z_0,...,z_m)\in I(\Ss_\g(e))\ \iff\ z_i\ne z_j \text{ for } i\ne j.
$$
\end{Theorem}
\begin{proof}
Note that in all these cases, if $e+F$ is semisimple then it is in fact \emph{regular}. This follows from the more general fact --- if $e$ is distinguished, and $e+F$ is semisimple, then $e+F$ is regular semisimple, see \cite{springer}*{9.5}. It follows that
$$
\Ss_\g(e)=\setof{F\in\g_{-d}}{p_r(e+F)\ne0},
$$
where $r$ is the rank of $\g$ and $p_r$ is the lowest nonzero coefficient (at degree $r$) of the characteristic polynomial of $\ad(e+F)$.

Obviously for irreducible nilpotent elements $e$ with $\dim\g_{-d}=1$, $e+F$ is semisimple if and only if $F\in\g_{-d}$ is nonzero, see Proposition \ref{propcon}.

When $\dim\g_{-d}=2$, there are exactly three distinct one-dimensional subspaces in $\g_{-d}$ such that $e+F$ is semisimple if and only if $F$ does not lie in any of those subspaces. We show this by a case-wise inspection of the four cases with $\dim\g_{-d}=2$ from Table \ref{irreds}.

\

\noindent{\bf Case 4$_k$: $\g=\so(4k+4)$, nilpotent element $e$ with partition $(2k+3,2k+1)$.}

The standard representation has a basis $x_{-k-1},x_{-k},...,x_{-1}$, $x_0$, $x_1,...,x_k$, $x_{k+1}$, $y_{-k},...,y_{-1}$, $y_0$, $y_1,...,y_k$, with $e$ acting by
\begin{equation}\label{basis4k}
\setlength\arraycolsep{2pt}
\begin{array}{lllllllllllllllllll}
x_{-k-1}&\mapsto&x_{-k}&\mapsto&\cdots&\mapsto&x_{-1}&\mapsto&x_0&\mapsto&x_1&\mapsto&\cdots&\mapsto&x_k&\mapsto&x_{k+1}&\mapsto&0\\
      &       &y_{-k}&\mapsto&\cdots&\mapsto&y_{-1}&\mapsto&y_0&\mapsto&y_1&\mapsto&\cdots&\mapsto&y_k&\mapsto&0.
\end{array}
\end{equation}
In this case $\g_{-d}$ has a basis $(F_1,F_2)$ such that
$$
F_1(y_k)=-x_{-k-1},\ F_1(x_{k+1})=y_{-k},\ F_2(x_k)=x_{-k-1},\ F_2(x_{k+1})=x_{-k},\
$$
and all other actions of $F_1$, $F_2$ are zero. Pictorially,
\begin{center}
\resizebox{!}{44ex}{
\begin{tikzpicture}[>=stealth,->,xscale=.05]
\useasboundingbox (-5,-5) rectangle (5,5);
\node (x0) at (0,0) {$x_0$};
\node (x1) [right=of x0] {$x_1$};
\node (xp) [right=of x1] {$\cdots$};
\node (xkm1) [right=of xp] {$x_k$};
\node (xk) [right=of xkm1] {$x_{k+1}$};
\node (xm1) [left=of x0] {$x_{-1}$};
\node (xm) [left=of xm1] {$\cdots$};
\node (x1mk) [left=of xm] {$x_{-k}$};
\node (xmk) [left=of x1mk] {$x_{-k-1}$};
\node (y0) [below=10em of x0] {$y_0$};
\node (y1) [below=10em of x1] {$y_1$};
\node (yp) [right=of y1] {$\cdots$};
\node (ykm1) [below=10em of xkm1] {$y_k$};
\node (ym1) [below=10em of xm1] {$y_{-1}$};
\node (ym) [left=of ym1] {$\cdots$};
\node (y1mk) [below=10em of x1mk] {$y_{-k}$};
\path[decoration=arrows, decorate] (xmk) -- node [above] {$e$} (x1mk) -- node [above] {$e$} (xm) -- node [above] {$e$} (xm1) -- node [above] {$e$} (x0) -- node [above] {$e$} (x1) -- node [above] {$e$} (xp) -- node [above] {$e$} (xkm1) -- node [above] {$e$} (xk);
\path[decoration=arrows, decorate] (y1mk) -- node [below] {$e$} (ym) -- node [below] {$e$} (ym1) -- node [below] {$e$} (y0) -- node [below] {$e$} (y1) -- node [below] {$e$} (yp) -- node [below] {$e$} (ykm1);
\draw[green,thick] (xkm1) .. controls +(up:10em) and +(up:20em) .. node[above,green] {$F_2$} (xmk);
\draw[green,thick] (xk) .. controls +(up:20em) and +(up:10em) .. node[above,green] {$F_2$} (x1mk);
\draw[red,thick] (xk) .. controls +(right:400em) and +(left:400em) .. node[above,red] {$F_1$} (y1mk);
\draw[red,thick] (ykm1) .. controls +(right:400em) and +(left:400em) .. node[above,red] {$-F_1$} (xmk);
\end{tikzpicture}
}
\end{center}

So $F=\lambda_1F_1+\lambda_2F_2$ acts via
\begin{align*}
x_k&\mapsto\lambda_2x_{-k-1},\\
y_k&\mapsto-\lambda_1x_{-k-1},\\
x_{k+1}&\mapsto\lambda_1y_{-k}+\lambda_2x_{-k},
\end{align*}
mapping all other $x_i$, $y_j$ to 0. Thus $c:=e+F$ acts as follows:
$$
\setlength\arraycolsep{1pt}
\begin{array}{lllllllll}
x_{-k-1}&\mapsto&x_{-k}&\mapsto&\cdots&\mapsto&x_k&\mapsto&x_{k+1}+\lambda_2x_{-k-1}\\
        &       &y_{-k}&\mapsto&\cdots&\mapsto&y_k&\mapsto&-\lambda_1x_{-k-1},\\
x_{k+1} &\mapsto&\lambda_1y_{-k}+\lambda_2x_{-k}&\mapsto&\cdots&\mapsto&\lambda_1y_k+\lambda_2x_k&\mapsto&\lambda_2x_{k+1}+(\lambda_2^2-\lambda_1^2)x_{-k-1}.
\end{array}
$$
\begin{align*}
c(y_{-k})&=y_{-k+1},\\
c^2(y_{-k})=c(y_{-k+1})&=y_{-k+2},\\
...\\
c^{2k}(y_{-k})=c(y_{k-1})&=y_k,\\
c^{2k+1}(y_{-k})=c(y_k)&=-\lambda_1x_{-k-1},\\
c^{2k+2}(y_{-k})=c(-\lambda_1x_{-k-1})&=-\lambda_1x_{-k},\\
c^{2k+3}(y_{-k})=c(-\lambda_1x_{-k})&=-\lambda_1x_{-k+1},\\
...\\
c^{4k+2}(y_{-k})=c(-\lambda_1x_{k-1})&=-\lambda_1x_k,\\
c^{4k+3}(y_{-k})=c(-\lambda_1x_k)&=-\lambda_1x_{k+1}-\lambda_1\lambda_2x_{-k-1},\\
\intertext{and}
c^{4k+4}(y_{-k})=c(-\lambda_1x_{k+1}-\lambda_1\lambda_2x_{-k-1})&=-\lambda_1^2y_{-k}-2\lambda_1\lambda_2x_{-k}\\
&=-\lambda_1^2y_{-k}-2\lambda_2c^{2k+2}(y_{-k}).
\end{align*}
In particular, since $c^{2k+1}(y_{-k})=-\lambda_1x_{-k-1}$ it follows that $c$ fails to be semisimple if $\lambda_1=0$. Whereas if $\lambda_1\ne0$, then $y_{-k},c(y_{-k}),c^2(y_{-k}),...,c^{4k+3}(y_{-k})$ form a basis of the standard representation, so that the action of $c$ on it can be realized as multiplication by $t$ on $\mathbb C[t]/(t^{4k+4}+2\lambda_2t^{2k+2}+\lambda_1^2)$. Discriminant of $t^{4k+4}+2\lambda_2t^{2k+2}+\lambda_1^2$ being a scalar multiple of $\lambda_1^{4k+2}(\lambda_1^2-\lambda_2^2)^{2k+2}$, we see that semisimplicity of $c$ can additionally fail only when $\lambda_2=\pm\lambda_1$. In this case it indeed fails since then $c^{2k+2}\pm\lambda_1\,\textrm{identity}$ becomes a nontrivial nilsquare operator.

So, semisimplicity of $e+F$ is equivalent to the conjunction of $\lambda_1\ne0$ and $\lambda_2\ne\pm\lambda_1$. Thus in this case the statement of the Theorem is ensured with the parametrization $\lambda_1=z_0-z_2$, $\lambda_2=z_0-2z_1+z_2$.

\

\noindent{\bf Case $\g=$F$_4$, nilpotent element $e$ with label F$_4(a_2)$}

Take the representative of this orbit
$$
e:=e_{1100}+e_{0011}+e_{0110}+e_{0210}
$$
(here $e_{ijkl}$ stands for the root vector of the root that is the linear combination of simple roots with coefficients $i$, $j$, $k$, $l$, where the numbering of simple roots is \raisebox{-.2ex}{\resizebox{!}{2ex}{\begin{tikzpicture}[>=stealth]\node (1) at (1.5,0) {1};\node (2) at (3,0) {2};\node (3) at (4.5,0) {3};\node (4) at (6,0) {4};\draw[thick] (1) -- (2); \draw[double,->,line width=.3mm] (2) -- (3); \draw[thick] (3) -- (4);\end{tikzpicture}}}).

Then $p_4(e+x_1f_{2431}+x_2f_{2432})$ is a scalar multiple of $x_1^2x_2^4(x_1+x_2)^2$, so the element $e+x_1f_{2431}+x_2f_{2432}$ is regular semisimple if and only if neither of the equalities $x_1=0$, $x_2=0$ or $x_1+x_2=0$ hold.

Obviously in this case the theorem holds true with $x_1=z_0-z_1$, $x_2=z_1-z_2$.

\

\noindent{\bf Case $\g=$E$_8$, nilpotent element $e$ with label E$_8(a_5)$.}
We take
$$
e:=
e_{\foreight{.5}00010000}+
e_{\foreight{.5}00000010}+
e_{\foreight{.5}10100000}+
e_{\foreight{.5}00110000}+
e_{\foreight{.5}01110000}+
e_{\foreight{.5}00111000}+
e_{\foreight{.5}01011100}+
e_{\foreight{.5}00001111}.
$$
The space $\g_{-d}$ has a basis consisting of negative root vectors
$$
F_1=f_{\foreight{.5}23465431},\ F_2=f_{\foreight{.5}23465432}.
$$
Then $p_8(e+x_1F_1+x_2F_2)$ is a scalar multiple of $x_1^8x_2^6(x_1+x_2)^6$, so that the theorem holds with the same parametrization as for the F$_4$ case above.

\

\noindent{\bf Case $\g=$E$_8$, nilpotent element $e$ with label E$_8(a_6)$.}
Here we take
$$
e:=
e_{\foreight{.5}00011000}+
e_{\foreight{.5}00000110}+
e_{\foreight{.5}00000011}+
e_{\foreight{.5}10110000}+
e_{\foreight{.5}01011000}+
e_{\foreight{.5}00011100}+
e_{\foreight{.5}00001110}+
e_{\foreight{.5}01111100}.
$$
The negative root vector basis $(F_1,F_2)$ here is the same as for E$_8(a_5)$, and $p_8(e+x_1F_1+x_2F_2)$ is a scalar multiple of $x_1^8(x_1-x_2)^8x_2^8$,
so that the theorem in this case is proved with the parametrization $x_1=z_0-z_1$, $x_2=z_2-z_1$.

\

There is only one case with $\dim\g_{-d}=3$: {\bf nilpotent element with label E$_7(a_5)$ in E$_7$.}

Take the representative
$$
e:=e_{\foreseven{.1}{.5}0000001}+e_{\foreseven{.1}{.5}1011000}+e_{\foreseven{.1}{.5}0111000}+e_{\foreseven{.1}{.5}0101100}+e_{\foreseven{.1}{.5}0011100}+e_{\foreseven{.1}{.5}0001110}+e_{\foreseven{.1}{.5}1111110}.
$$

Let
$$
c(x_1,x_2,x_3):=e+x_1f_{\foreseven{.1}{.5}1224321}+x_2f_{\foreseven{.1}{.5}1234321}+x_3f_{\foreseven{.1}{.5}2234321},
$$
then $p_7(c(x_1,x_2,x_3))$ is a scalar multiple of
$$
(x_1-x_3)^4(x_1^2+x_1 x_3+x_3^2)^4 (x_1-x_2+x_3)^3 (x_1^2+x_1 x_2-x_1 x_3+x_2^2+x_2 x_3+x_3^2)^3.
$$
Denoting by $\omega$ the primitive third root of unity, we have
$$
x_1^2+x_1 x_3+x_3^2=(x_1-\omega x_3)(x_1-\bar\omega x_3)
$$
and
$$
x_1^2+x_1 x_2-x_1 x_3+x_2^2+x_2 x_3+x_3^2=(x_2-\omega x_1-\bar\omega x_3)(x_2-\bar\omega x_1-\omega x_3),
$$
so that semisimplicity of $c(x_1,x_2,x_3)$ fails along the following subset of the projective plane:

\begin{center}
\resizebox{.75\textwidth}{!}{%
\begin{tikzpicture}[scale=2]
\coordinate (A) at (1,0);
\coordinate (B) at ({cos(120)},{sin(120)});
\coordinate (C) at ({cos(240)},{sin(240)});
\coordinate (O) at (0,0);
\coordinate (a) at (-1/2,0);
\coordinate (b) at ({-cos(120)/2},{-sin(120)/2});
\coordinate (c) at ({-cos(240)/2},{-sin(240)/2});
\draw [add= 1 and 1, green, ultra thick] (A) to (B);
\draw [add= 1 and 1, green, ultra thick] (A) to (C);
\draw [add= 1 and 1, green, ultra thick] (B) to (C);
\draw [add= 1 and 1, blue, ultra thick] (A) to (a);
\draw [add= 1 and 1, blue, ultra thick] (B) to (b);
\draw [add= 1 and 1, blue, ultra thick] (C) to (c);
\node[circle,fill,red] at (A) {};
\node[circle,fill,red] at (B) {};
\node[circle,fill,red] at (C) {};
\node[circle,fill,brown] at (a) {};
\node[circle,fill,brown] at (b) {};
\node[circle,fill,brown] at (c) {};
\node[circle,fill] at (O) {};
\node at (-2.5,0) {$x_1=x_3$};
\node at (2.5,-1) {$x_2=\bar\omega x_1+\omega x_3$};
\node at (-2.5,-1.5) {$x_2=\omega x_1+\bar\omega x_3$};
\node at (-1/2,2.75) {$x_2=x_1+x_3$};
\node at (1,2) {$x_1=\bar\omega x_3$};
\node at (1,-2) {$x_1=\omega x_3$};
\node at (1,.22) {${}_{[1:-1:1]}$};
\node at (-.1,.95) {${}_{[\bar\omega:-1:\omega]}$};
\node at (-.1,-1) {${}_{[\omega:-1:\bar\omega]}$};
\node at (.3,.1) {${}_{[0:1:0]}$};
\node at (-.8,.15) {${}_{[1:2:1]}$};
\node at (.6,-.45) {${}_{[\bar\omega:2:\omega]}$};
\node at (.6,.45) {${}_{[\omega:2:\bar\omega]}$};
\end{tikzpicture}%
}
\end{center}

For this case we can ensure the theorem with
$$
x_1=\frac{\omega z_1+\omega^2z_2+z_3}3,\ x_2=z_0-\frac{z_1+z_2+z_3}3,\ x_3=\frac{\omega^2z_1+\omega z_2+z_3}3.
$$

\

Finally, for $\dim(\g_{-d})=4$ there is also only one case: {\bf nilpotent orbit labeled by E$_8(a_7)$ in E$_8$.}

The theorem in this case has been inspired by an answer that Noam Elkies gave to a question on mathoverflow concerning the configuration of hyperplanes that appears in this case --- see \cite{MOElkies}.

We take
$$
e:=
e_{\foreight{.5}00001000}+
e_{\foreight{.5}00111100}+
e_{\foreight{.5}10111100}+
e_{\foreight{.5}01011110}+
e_{\foreight{.5}00011111}+
e_{\foreight{.5}11111100}+
e_{\foreight{.5}01111110}+
e_{\foreight{.5}11221000}.
$$
The root vector basis of $\g_{-d}$ consists of negative root vectors
$$
F_1:=f_{\foreight{.5}23465321},\ F_2:=f_{\foreight{.5}23465421},\ F_3:=f_{\foreight{.5}23465431},\ F_4:=f_{\foreight{.5}23465432}.
$$
Here $p_8(e+x_1F_1+x_2F_2+x_3F_3+x_4F_4)$ is a scalar multiple of the 24th power of
$$
x_1x_3x_4(x_3+x_4)(x_1^2+x_2^2)(x_1^2+(x_2+x_3)^2)(x_1^2+(x_2-x_4)^2).
$$
Replacing $x_1$ with $\sqrt{-1}x_1$, we find that the singular set consists of ten $3$-di\-men\-si\-o\-nal subspaces of $\g_{-d}$, given in the root vector basis by the equations
\begin{align*}
&x_1-x_2-x_3=0, x_1+x_2+x_3=0, x_1-x_2+x_4=0, x_1+x_2-x_4=0,\\
&x_1=0, x_3=0, x_4=0, x_1+x_2=0, x_1-x_2=0, x_3+x_4=0.
\end{align*}
All possible intersections of these subspaces produce twenty five 2-dimensional subspaces and fifteen $1$-dimensional subspaces. Each $3$-dimensional subspace contains six of these $2$-dimensional subspaces and seven of these $1$-dimensional subspaces. Each of these $2$-dimensional subspaces contains three of the $1$-dimensional subspaces. Ten of the $1$-dimensional subspaces lie in four of the $2$-dimensional and in four of the $3$-dimensional subspaces each, while five of the $1$-dimensional subspaces lie in seven of the $2$-dimensional and in six of the $3$-dimensional subspaces each. Finally fifteen of the $2$-dimensional subspaces lie in two of the $3$-dimensional ones and ten of the $2$-dimensional subspaces lie in three of the $3$-dimensional ones.

The parametrization (found by Noam Elkies in \cite{MOElkies}) in this case is
$$
x_1=z_0 - z_1,\quad x_2=z_0 + z_1 - 2z_2,\quad x_3=2(z_2 - z_3),\quad x_4=2(z_4 - z_2).
$$
This parametrization in particular shows that the whole configuration can be described through its projectivization as the barycentric subdivision of a tetrahedron:

\begin{center}
\includegraphics[scale=.18,trim=0 20 0 30,clip]{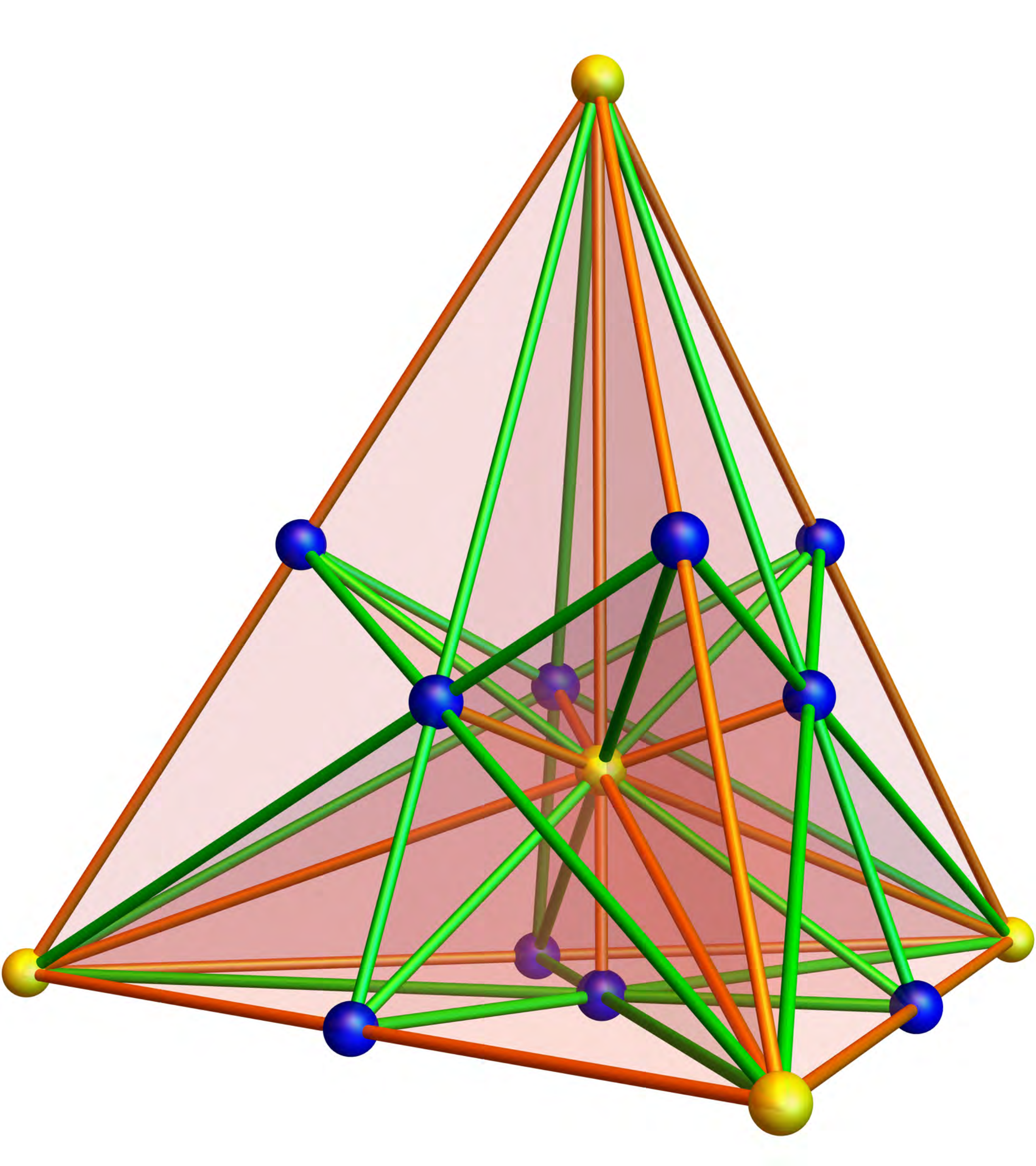}
\end{center}

The above fifteen $1$-dimensional subspaces correspond to its vertices ($4$), ba\-ry\-cen\-ters of edges ($6$), barycenters of faces ($4$) and the barycenter of the tetrahedron ($1$), twenty five $2$-dimensional subspaces correspond to edges ($6$), lines joining a vertex with the barycenter of some face ($4\times 4$) and lines joining barycenters of opposite edges ($3$), and ten $3$-dimensional subspaces of the configuration correspond to faces ($4$) and planes through an edge and the barycenter of the tetrahedron ($6$).
\end{proof}

\section{Non-irreducible nilpotent elements of semisimple type}\label{sectionred}

As shown in \cite{ekv}*{Theorem 3.14}, for each nilpotent element $e$ of semisimple type there is a reducing subalgebra for $e$ where it is of \emph{regular} semisimple type. We will, in fact, for each such $e$ exhibit a reducing subalgebra where it is irreducible (hence regular).

In most cases, these reducing subalgebras are as follows.

\begin{Definition}
For a nilpotent element $e$, let $\la(e)$ denote the semisimple part of the centralizer of a Cartan subalgebra of the centralizer $\z(\sa)$ of the $\sla(2)$-triple $\sa$ for $e$.
\end{Definition}

The subalgebra $\la(e)$ is the derived subalgebra of a minimal Levi subalgebra of $\g$ containing $\sa$, and $e$ is distinguished in it, so that $\sa$ has zero centralizer in $\la(e)$. It turns out, by looking at Tables \ref{reductions}, \ref{reductionse6}, \ref{reductionse7}, \ref{reductionse8}, \ref{reductionsf4} that for most of nilpotent elements $e$ in $\g$ of semisimple type, $\la(e)$ is a reducing subalgebra for $e$. The exceptions in classical $\g$, when $\la(e)$ is not a reducing subalgebra, are the following (see \cite{ekv}, before Section 5):
\begin{itemize}
\item[(a)] nilpotent elements with partition $[3,1^{(2k)}]$ in $\so(2k+3)$ for $k\ge1$, $\rank=2$,
\item[(b)] nilpotent elements with partition $[(2k)^{(n)}]$ in $\spa(2kn)$ for $n>1$, $k\ge1$, $\rank=n$.
\end{itemize}

In case (a), the algebra $\la(e)$ has type A$_1$, with $\la(e)_{-d}$ of dimension 1, while $\g_{-d}$ has dimension $2k+1$ and $\rank e=2$. The centralizer of $\sa$ in $\g$ is $\so(2k)$ acting trivially on $\la(e)_{-d}$, so $\la(e)$ cannot be reducing.

In case (b), $\g_{-d}$ has dimension $n(n+1)/2$, with the centralizer $\so(n)$ of $\sa$ acting on $\g_{-d}$ as on the symmetric square of the standard representation, so that $e$ has rank $n$, while $\la(e)$ is $\sla(2k)^{\oplus j}$ for $n=2j$ and $\sla(2k)^{\oplus(j-1)}\oplus\spa(2k)$ for $n=2j-1$, with $e$ principal, hence of rank $j$ in $\la(e)$ in both cases.

There is only one nilpotent element $e$ in exceptional $\g$, when the algebra $\la(e)$ is not reducing, namely for $e$ with label $\tilde\rA_1$ in F$_4$, which has rank 2. Here the centralizer of $\sa$ is $\sla(4)$, and $\g_{-d}$ is the sum of a 6-dimensional irreducible $\sla(4)$-module and a 1-dimensional trivial module. Since $e$ has rank 1 in $\la(e)$, the latter cannot be a reducing subalgebra.

In these three cases, minimal reducing subalgebras are the ones generated by $e$ and an element $F\in\g_{-d}$ having closed orbit of smallest possible codimension (equal to the rank of the nilpotent). In case (a) and for $\tilde{\text A}_1$ in F$_4$ it is of type A$_1+$A$_1$, and in case (b) it is $\spa(2k)^{\oplus j}$.

In all remaining cases, $\la(e)$ is reducing, and $e$ is principal in $\la(e)$.

There are also several cases when, although $\la(e)$ is a reducing subalgebra, there is a still smaller reducing subalgebra inside it. Such subalgebra is generated by $e$ and an element $F\in\la(e)_{-d}$ as above --- that is, an element having closed orbit of smallest possible codimension. In all these cases it turns out that $e$ is irreducible in this subalgebra, i.~e. it gives one of the cases from Table \ref{irreds}. (We have only a computer proof of this.) It then follows that this is a \emph{minimal} reducing subalgebra.

Thus in Tables \ref{reductions}, \ref{reductionse6}, \ref{reductionse7}, \ref{reductionse8}, \ref{reductionsf4} all algebras in the column ``$\qa_{\min}$'' are minimal reducing subalgebras, and have the property that they are generated by $e$ and $F\in\g_{-d}$, having closed orbit of minimal codimension.

In Tables \ref{reductions}, \ref{reductionse6}, \ref{reductionse7}, \ref{reductionse8}, \ref{reductionsf4} we list all nilpotent orbits $G(e)$ of semisimple type, except for the irreducible ones, in all simple Lie algebras (the irreducible ones are listed in Table \ref{irreds}). In the first column the nilpotent elements are given by the corresponding partitions in the classical Lie algebras (notation $k^{(s)}$ means that the part $k$ is repeated $s$ times), and by the type of $\la(e)$ and by the weighted Dynkin diagram in the exceptional Lie algebras. In the second and third columns the depth and rank are given. In the fourth column the image of $Z(\sa)$ in End$(\g_{-d})$ is given. It is computed using the $\z(\sa)$ listed in \cite{ekv} and the results of \cites{ale,CM}. Actions of $\z(\sa)$ on $\g_{-d}$ are computed using the GAP command $\mathtt{LeftAlgebraModule}$ \cite{gap} which finds the module structure. For the torus part of $\z(\sa)$ one finds eigenvectors and eigenvalues of its action on $\g_{-d}$. Next, the command $\mathtt{DirectSumDecomposition}$ in \cite{gap} decomposes $\g_{-d}$ as a module over the semisimple part of $\z(\sa)$ into irreducible components. For all $\g$ of exceptional type dimensions of these irreducible components suffice to determine the structure of these irreducible components up to isomorphism. In the fifth column the minimal reducing subalgebras are given by their number in Table \ref{irreds} (recall that $e$ is irreducible in its minimal reducing subalgebra).

Concerning notation --- ``st'' denotes the standard representations, ``ad'' the adjoint representations, {\bf1} the trivial 1-dimensional representations, {\bf7} and {\bf26} the non-trivial irreducible representations of minimal dimension of G$_2$ and F$_4$ respectively, $\sigma_n$ is the nontrivial irreducible $n-1$-dimensional representation of the symmetric group S$_n$ ($n\ge2$), $\sigma_n\oplus{\bf1}$ being its permutation representation, and $-^{\oplus n}$ is the direct sum with itself $n$ times. In all cases $k\ge1$, $n\ge1$, $q\ge0$.

We also list $\qa_{\max}$, which is a subalgebra of $\g$, generated by $\g_{-d}$ and $e$ (it is a reducing subalgebra by the results of \cite{ekv}). Types of $\qa_{\min}$ and $\qa_{\max}$ are determined using the GAP command $\mathtt{SemiSimpleType}$ \cite{gap}. Finally, in the last column we list the algebras $(\g_{-d},*)$ (their notation is explained in Section \ref{algsection}). They are defined by \eqref{multab} and in Section {\bf5C}. As in Table \ref{irreds}, $1$ stands for the $1$-dimensional algebra with non-zero (resp.~zero) multiplication if $\frac d2$ is odd (resp.~even).

\begin{center}
\begin{table}[H]\renewcommand\thetable{2ABCD}\caption{Non-irreducible nilpotent elements of semisimple type in A, B, C, D\label{reductions}}

\

\setlength\tabcolsep{1pt}
\resizebox{\textwidth}{!}{
\begin{tabular}{l|r|r|r|l|l|rr}
nilpotent $e$&depth&rank&$Z(\sa)|\g_{-d}$&$\qa_{\min}$&$\qa_{\max}\ni e$&$(\g_{-d},*)$\\
\hline
\multicolumn{4}{l|}{$\sla$}\\
\hline
$[(2k)^{(n)},1^{(q)}]$&$4k-2$&$n$&ad$_{\sla(n)}\oplus\bf1$&$2_k^{\oplus n}$&$\sla(2kn)\ni[(2k)^{(n)}]$&$J_n(A)$\bigstrut[t]\\
$[(2k+1)^{(n)},1^{(q)}]$&$4k$&$n$&ad$_{\sla(n)}\oplus\bf1$&$1_k^{\oplus n}$&$\sla(2kn+n)\ni[(2k+1)^{(n)}]$&A$_{n-1}\oplus{\bf1}$\bigstrut\\
\hline
\multicolumn{4}{l|}{$\spa$}\\
\hline
$[(2k)^{(n)},1^{(2q)}]$&$4k-2$&$n$&S$^2($st$_{\so(n)})$&$2_k^{\oplus n}$&$\spa(2kn)\ni[(2k)^{(n)}]$&$J_n(C)$\bigstrut[t]\\
$[(2k+1)^{(2n)},1^{(2q)}]$&$4k$&$n$&ad$_{\spa(2n)}$&$1_k^{\oplus n}$&$\spa(4nk+2n)\ni[(2k+1)^{(2n)}]$&C$_n$\bigstrut\\
\hline
\multicolumn{4}{l|}{$\so$}\\
\hline
$[3,1^{(q)}]$, $q\ne1$&$2$&2&st$_{\so(q)}\oplus\bf1$&$2_1^{\oplus2}$&$\g$&$J_q(BD)$\bigstrut[t]\\
$[7,1^{(q)}]$&$10$&1&{\bf1}&$5$&$=\qa_{\min}$&{\bf1}\\
$[2k+3,1^{(q)}]$, $k\ne2$&$4k+2$&1&{\bf1}&$3_{k+1}$&$=\qa_{\min}$&{\bf1}\\
$[2k+3,2k+1,1^{(q)}]$&$4k+2$&2&${\bf1}\oplus{\bf1}$&$4_k$&$=\qa_{\min}$&$\CA_{-k}(2)$\\
$[(2k)^{(2)},1^{(q)}]$&$4k-2$&$1$&$\bf1$&$2_k$&$=\qa_{\min}$&{\bf1}\\
$[(2k)^{(2n)},1^{(q)}]$, $n>1$&$4k-2$&$n$&$\Lambda^2(\text{st}_{\spa(2n)})$&$2_k^{\oplus n}$&$\so(4kn)\ni[(2k)^{(2n)}]$&$J_{2n}(D)$\\
$[(2k+1)^{(2)},1^{(q)}]$&$4k$&$1$&$\bf1$&$1_k$&$=\qa_{\min}$&{\bf1}\\
$[(2k+1)^{(2n)},1^{(q)}]$, $n>1$&$4k$&$n$&ad$_{\so(2n)}$&$1_k^{\oplus n}$&$\so(4nk+2n)\ni[(2k+1)^{(2n)}]$&D$_n$
\end{tabular}
}
\end{table}

\begin{table}[H]\renewcommand\thetable{2FG}\caption{Non-irreducible nilpotent elements of semisimple type in F$_4$ and G$_2$\label{reductionsf4}}

\

\resizebox{\textwidth}{!}{
\begin{tabular}{ll|r|r|r|l|l|rr}
&nilpotent $e$&depth&rank&$Z(\sa)|\g_{-d}$&$\qa_{\min}$&$\qa_{\max}\ni e$&$(\g_{-d},*)$\\
\hline
\multicolumn{5}{l|}{F$_4$}\\
\hline
&A$_1$\hfill\FDynkin{0,1,0,0}&2&1&{\bf1}&$2_1$&$=\qa_{\min}$&{\bf1}\bigstrut[t]\\
&$\tilde{\mathrm A}_1$\textsuperscript{$1)$}\hfill\FDynkin{1,0,0,0}
&2&2&st$_{\so(6)\rtimes\mathrm S_2}\oplus\bf1$&$2_1^{\oplus2}$&$\so(9)\ni[3,1^{(6)}]$&$J_6(BD)$\\
&A$_2$\hfill\FDynkin{0,2,0,0}&4&1&$\sigma_2$&$1_1$&$=\qa_{\min}$&{\bf1}\\
&$\tilde{\text A}_2$\hfill\FDynkin{2,0,0,0}&4&1&${\bf7}$&$1_1$&$\g$&$M$\\
&B$_2$\hfill\FDynkin{1,2,0,0}&6&1&{\bf1}&$3_2$&$=\qa_{\min}$&{\bf1}\\
&F$_4(a_3)$\hfill\FDynkin{0,0,0,2}&6&2&$\sigma_3$&$4_1$&$=\qa_{\min}$&$\CA_{-1}(2)$\\
&B$_3$\hfill\FDynkin{0,2,0,2}&10&1&{\bf1}&$5$&$=\qa_{\min}$&{\bf1}\\
&C$_3$\hfill\FDynkin{2,1,1,0}&10&1&{\bf1}&$2_3$&$=\qa_{\min}$&{\bf1}\\
&F$_4(a_1)$\hfill\FDynkin{2,2,0,2}&14&1&{\bf1}&$3_4$&$=\qa_{\min}$&{\bf1}\\
\hline
\multicolumn{5}{l|}{G$_2$}\\
\hline
&A$_1$\hfill\GDynkin{0,1}&2&1&${\bf1}$&$2_1$&$=\qa_{\min}$&{\bf1}\bigstrut[t]\\
&G$_2(a_1)$\hfill\GDynkin{0,2}&4&1&${\bf1}$&$1_1$&$=\qa_{\min}$&{\bf1}
\end{tabular}
}
\end{table}
\end{center}

\footnotetext[1]{\!$^{)}$ Here the action of S$_2$ on the standard representation of $\so(6)$ is the one which induces the non-trivial diagram automorphism of $\so(6)$}

\begin{table}[H]\renewcommand\thetable{2E6}\caption{Non-irreducible nilpotent elements of semisimple type in E$_6$\label{reductionse6}}

\

\resizebox{\textwidth}{!}{
\begin{tabular}{l|r|r|r|l|l|rr}
nilpotent $e$&depth&rank&$Z(\sa)|\g_{-d}$&$\qa_{\min}$&$\qa_{\max}\ni e$&$(\g_{-d},*)$\\
\hline
A$_1$\hfill\EDynkin{1,0,0,0,0,0}{1}&2&1&{\bf1}&$2_1$&$=\qa_{\min}$&{\bf1}\bigstrut[t]\\
2A$_1$\hfill\EDynkin{0,1,0,0,0,1}{1}&2&2&st$_{\so(7)}\oplus\bf1$&$2_1^{\oplus2}$&$\so(10)\ni[3,1^{(7)}]$&$J_7(BD)$\\
A$_2$\hfill\EDynkin{2,0,0,0,0,0}{1}&4&1&$\sigma_2$&$1_1$&$=\qa_{\min}$&{\bf1}\\
2A$_2$\hfill\EDynkin{0,2,0,0,0,2}{1}&4&2&${\bf7}\oplus\bf1$&$1_1^{\oplus2}$&$\g$&$M\oplus{\bf1}$\\
A$_3$\hfill\EDynkin{2,1,0,0,0,1}{1}&6&1&${\bf1}$&$3_2$&$=\qa_{\min}$&{\bf1}\\
D$_4(a_1)$\hfill\EDynkin{0,0,0,2,0,0}{1}&6&2&$\sigma_3$&$4_1$&$=\qa_{\min}$&$\CA_{-1}(2)$\\
A$_4$\hfill\EDynkin{2,2,0,0,0,2}{1}&8&1&${\bf1}$&$1_2$&$=\qa_{\min}$&{\bf1}\\
D$_4$\hfill\EDynkin{2,0,0,2,0,0}{1}&10&1&${\bf1}$&$5$&$=\qa_{\min}$&{\bf1}\\
A$_5$\hfill\EDynkin{1,2,1,0,1,2}{1}&10&1&${\bf1}$&$2_3$&$=\qa_{\min}$&{\bf1}\\
E$_6(a_3)$\hfill\EDynkin{0,2,0,2,0,2}{1}&10&2&$\sigma_2\oplus{\bf1}$&$7$&$=\qa_{\min}$&$\CA_{-\frac13}(2)$\\
D$_5$\hfill\EDynkin{2,2,0,2,0,2}{1}&14&1&${\bf1}$&$3_4$&$=\qa_{\min}$&{\bf1}\\
E$_6$\hfill\EDynkin{2,2,2,2,2,2}{1}&22&1&${\bf1}$&$6$&$=\qa_{\min}$&{\bf1}
\end{tabular}
}
\end{table}

\begin{table}[H]\renewcommand\thetable{2E7}\caption{Non-irreducible nilpotent elements of semisimple type in E$_7$\label{reductionse7}}

\

\resizebox{\textwidth}{!}{
\begin{tabular}{l|r|r|r|l|l|rr}
nilpotent $e$&depth&rank&$Z(\sa)|\g_{-d}$&$\qa_{\min}$&$\qa_{\max}\ni e$&$(\g_{-d},*)$\\
\hline
A$_1$\hfill\EDynkin{0,1,0,0,0,0,0}{1}&2&1&{\bf1}&$2_1$&$=\qa_{\min}$&{\bf1}\bigstrut[t]\\
2A$_1$\hfill\EDynkin{0,0,0,0,0,1,0}{1}&2&2&st$_{\so(9)}\oplus\bf1$&$2_1^{\oplus2}$&$\so(12)\ni[3,1^{(9)}]$&$J_9(BD)$\\
$[3\rA_1]''$\hfill\EDynkin{0,0,0,0,0,0,2}{1}&2&3&${\bf26}\oplus\bf1$&$2_1^{\oplus3}$&$\g$&$J(E)$\\
A$_2$\hfill\EDynkin{0,2,0,0,0,0,0}{1}&4&1&$\sigma_2$&$1_1$&$=\qa_{\min}$&{\bf1}\\
2A$_2$\hfill\EDynkin{0,0,0,0,0,2,0}{1}&4&2&${\bf7}\oplus$st$_{\so(3)}$&$1_1^{\oplus2}$&$\g$&$M\oplus$A$_1$\\
A$_3$\hfill\EDynkin{0,2,0,0,0,1,0}{1}&6&1&{\bf1}&$3_2$&$=\qa_{\min}$&{\bf1}\\
D$_4(a_1)$\hfill\EDynkin{0,0,2,0,0,0,0}{1}&6&2&$\sigma_3$&$4_1$&$=\qa_{\min}$&$\CA_{-1}(2)$\\
A$_4$\hfill\EDynkin{0,2,0,0,0,2,0}{1}&8&1&$\sigma_2$&$1_2$&$=\qa_{\min}$&{\bf1}\\
D$_4$\hfill\EDynkin{0,2,2,0,0,0,0}{1}&10&1&{\bf1}&$5$&$=\qa_{\min}$&{\bf1}\\
$[\rA_5]''$\hfill\EDynkin{0,2,0,0,0,2,2}{1}&10&1&{\bf1}&$2_3$&$=\qa_{\min}$&{\bf1}\\
$[\rA_5]'$\hfill\EDynkin{0,1,0,1,0,2,0}{1}&10&1&{\bf1}&$2_3$&$=\qa_{\min}$&{\bf1}\\
D$_6(a_2)$\hfill\EDynkin{1,0,1,0,1,0,2}{1}&10&2&${\bf1}\oplus\bf1$&$4_2$&$=\qa_{\min}$&$\CA_{-2}(2)$\\
E$_6(a_3)$\hfill\EDynkin{0,0,2,0,0,2,0}{1}&10&2&$\sigma_2\oplus{\bf1}$&$7$&$=\qa_{\min}$&$\CA_{-\frac13}(2)$\\
A$_6$\hfill\EDynkin{0,0,0,2,0,2,0}{1}&12&1&st$_{\so(3)}$&$1_3$&$\g$&A$_1$\\
D$_5$\hfill\EDynkin{0,2,2,0,0,2,0}{1}&14&1&{\bf1}&$3_4$&$=\qa_{\min}$&{\bf1}\\
E$_6(a_1)$\hfill\EDynkin{0,2,0,2,0,2,0}{1}&16&1&$\sigma_2$&$8$&$=\qa_{\min}$&{\bf1}\\
D$_6$\hfill\EDynkin{1,2,1,0,1,2,2}{1}&18&1&{\bf1}&$3_5$&$=\qa_{\min}$&{\bf1}\\
E$_6$\hfill\EDynkin{0,2,2,2,0,2,0}{1}&22&1&{\bf1}&$6$&$=\qa_{\min}$&{\bf1}
\end{tabular}
}
\end{table}

\begin{table}[H]\renewcommand\thetable{2E8}\caption{Non-irreducible nilpotent elements of semisimple type in E$_8$\label{reductionse8}}

\

\begin{center}
\resizebox{.9\textwidth}{!}{
\begin{tabular}{l|r|r|r|l|l|rr}
nilpotent $e$&depth&rank&$Z(\sa)|\g_{-d}$&$\qa_{\min}$&$\qa_{\max}\ni e$&$(\g_{-d},*)$\\
\hline
A$_1$\hfill\EDynkin{0,0,0,0,0,0,0,1}{1}&2&1&{\bf1}&$2_1$&$=\qa_{\min}$&{\bf1}\bigstrut[t]\\
2A$_1$\hfill\EDynkin{0,1,0,0,0,0,0,0}{1}&2&2&st$_{\so(13)}\oplus\bf1$&$2_1^{\oplus2}$&$\so(16)\ni[3,1^{(13)}]$&$J_{13}(BD)$\\
A$_2$\hfill\EDynkin{0,0,0,0,0,0,0,2}{1}&4&1&$\sigma_2$&$1_1$&$=\qa_{\min}$&{\bf1}\\
2A$_2$\hfill\EDynkin{0,2,0,0,0,0,0,0}{1}&4&2&${\bf7}\otimes(\sigma_2\oplus{\bf1})$&$1_1^{\oplus2}$&$\g$&$M\oplus M$\\
A$_3$\hfill\EDynkin{0,1,0,0,0,0,0,2}{1}&6&1&{\bf1}&$3_2$&$=\qa_{\min}$&{\bf1}\\
D$_4(a_1)$\hfill\EDynkin{0,0,0,0,0,0,2,0}{1}&6&2&$\sigma_3$&$4_1$&$=\qa_{\min}$&$\CA_{-1}(2)$\\
A$_4$\hfill\EDynkin{0,2,0,0,0,0,0,2}{1}&8&1&$\sigma_2$&$1_2$&$=\qa_{\min}$&{\bf1}\\
D$_4$\hfill\EDynkin{0,0,0,0,0,0,2,2}{1}&10&1&{\bf1}&$5$&$=\qa_{\min}$&{\bf1}\\
A$_5$\hfill\EDynkin{0,2,0,0,0,1,0,1}{1}&10&1&{\bf1}&$2_3$&$=\qa_{\min}$&{\bf1}\\
E$_6(a_3)$\hfill\EDynkin{0,2,0,0,0,0,2,0}{1}&10&2&$\sigma_2\oplus{\bf1}$&$7$&$=\qa_{\min}$&$\CA_{-1}(2)$\\
D$_6(a_2)$\hfill\EDynkin{1,0,1,0,0,0,1,0}{1}&10&2&$\sigma_2\oplus{\bf1}$&$4_2$&$=\qa_{\min}$&$\CA_{-2}(2)$\\
E$_7(a_5)$\hfill\EDynkin{0,0,0,1,0,1,0,0}{1}&10&3&$\sigma_3\oplus{\bf1}$&$11$&$=\qa_{\min}$&$\CA_{-\frac13}(3)$\\
A$_6$\hfill\EDynkin{0,2,0,0,0,2,0,0}{1}&12&1&st$_{\so(3)}$&$1_3$&E$_7\ni$A$_6$&A$_1$\\
D$_5$\hfill\EDynkin{0,2,0,0,0,0,2,2}{1}&14&1&{\bf1}&$3_4$&$=\qa_{\min}$&{\bf1}\\
E$_6(a_1)$\hfill\EDynkin{0,2,0,0,0,2,0,2}{1}&16&1&$\sigma_2$&$8$&$=\qa_{\min}$&{\bf1}\\
D$_6$\hfill\EDynkin{1,2,1,0,0,0,1,2}{1}&18&1&{\bf1}&$3_5$&$=\qa_{\min}$&{\bf1}\\
E$_6$\hfill\EDynkin{0,2,0,0,0,2,2,2}{1}&22&1&{\bf1}&$6$&$=\qa_{\min}$&{\bf1}\\
D$_7$\hfill\EDynkin{1,2,1,0,1,1,0,1}{1}&22&1&{\bf1}&$3_6$&$=\qa_{\min}$&{\bf1}\\
E$_7(a_1)$\hfill\EDynkin{1,2,1,0,1,0,2,2}{1}&26&1&{\bf1}&$10$&$=\qa_{\min}$&{\bf1}\\
E$_7$\hfill\EDynkin{1,2,1,0,1,2,2,2}{1}&34&1&{\bf1}&$9$&$=\qa_{\min}$&{\bf1}
\end{tabular}
}
\end{center}
\end{table}

\section{Reformulation in terms of algebra structure in $\g_{-d}$}\label{algsection}

\noindent{\bf 5A.}
Let $e$ be a nilpotent element in $\g$ of even depth $d$. Consider the binary operation
$$
[(\ad e)^i(x),(\ad e)^j(y)],\qquad x,y\in\g_{-d}.
$$
Since with respect to the grading \eqref{grading} defined by $e$, $e$ itself is homogeneous of degree $2$, clearly when $x$ and $y$ are both homogeneous of degree $-2(i+j)$, the result will be homogeneous of the same degree. Moreover for $i>0$ we have
\begin{multline*}
[(\ad e)^{i-1}(x),(\ad e)^{j+1}(y)]=[(\ad e)^{i-1}(x),[e,(\ad e)^j(y)]]\\=-[(\ad e)^j(y),[(\ad e)^{i-1}(x),e]]-[e,[(\ad e)^j(y),(\ad e)^{i-1}(x)]]\\
=[[(\ad e)^{i-1}(x),e],(\ad e)^j(y)]-[e,[(\ad e)^j(y),(\ad e)^{i-1}(x)]]\\=-[(\ad e)^i(x),(\ad e)^j(y)]-[e,[(\ad e)^j(y),(\ad e)^{i-1}(x)]].
\end{multline*}
Now, for $x,y\in\g_{-2(i+j)}$, the element $[(\ad e)^j(y),(\ad e)^{i-1}(x)]$ lies in $\g_{-2(i+j+1)}$, so that if $2(i+j)$ is equal to the depth $d$,
the latter element will be zero by dimension considerations. Hence we have, provided that $d$ is even,
\begin{equation}\label{plusminus}
[(\ad e)^{\frac d2}(x),y]=-[(\ad e)^{\frac d2-1}(x),(\ad e)(y)]=...=(-1)^{\frac d2}[x,(\ad e)^{\frac d2}(y)].
\end{equation}
It follows that all the operations that can be obtained in this way on $\g_{-d}$ differ only by sign. We will pick one of these and will always use the operation
\begin{equation}\label{operdef}
x*y:=[(\ad e)^{\frac d2}(x),y],\quad x,y\in\g_{-d}.
\end{equation}
It follows from \eqref{plusminus} that this operation is skew-commutative when $\frac d2$ is even and commutative when $\frac d2$ is odd
(for odd $d$ we do not get any operation on $\g_{-d}$).

Note that the $*$-algebra structure \eqref{operdef} is $Z(\sa)$-invariant. Note also that we have
\begin{Proposition}\label{symbil}
The symmetric bilinear form on $\g_{-d}$ given by
$$
(x,y)=\abr{(\ad e)^dx,y},
$$
where $\abr{\cdot,\cdot}$ is the Killing form, is non-degenerate and associative for the product \eqref{operdef}, provided that $d$ is even.
\end{Proposition}
\begin{proof}
Let us abbreviate the operator $\operatorname{ad}e$ to $E$, and $\operatorname{ad}f$ to $F$, where $\{e,f,h\}$ is the standard $\mathfrak{sl}(2)$-triple $\mathfrak s$. We have (by associativity of the Killing form)
\[
(x,y)=\abr{E^dx,y}=-\abr{E^{d-1}x,Ey}=...=(-1)^i\abr{E^{d-i}x,E^iy}=...=\abr{x,E^dy},
\]
and, by \eqref{plusminus},
\begin{multline*}
x*y=[E^{d/2}x,y]=-[E^{d/2-1}x,Ey]=...=(-1)^j[E^{d/2-j}x,E^jy]=...\\=(-1)^{d/2}[x,E^{d/2}y].
\end{multline*}
Thus to prove
\[
(x*y,z)=(x,y*z)
\]
means to prove
\[
\abr{E^d[E^{d/2}x,y],z}=\abr{E^dx,[E^{d/2}y,z]}.
\]

Let us transform the left hand side as
\[
\abr{E^d[E^{d/2}x,y],z}=\abr{[E^{d/2}x,y],E^dz}=-\abr{y,[E^{d/2}x,E^dz]},
\]
and the right hand side as
\[
\abr{E^dx,[E^{d/2}y,z]}=(-1)^{d/2}\abr{E^dx,[y,E^{d/2}z]}=-(-1)^{d/2}\abr{y,[E^dx,E^{d/2}z]}.
\]
We then see that it suffices (but in fact it is also easy to see that it is necessary) to prove
\[
[E^{d/2}x,E^dz]=(-1)^{d/2}[E^dx,E^{d/2}z]\qquad\text{for any $x,z\in\g_{-d}$.}
\]

Note that both $x$ and $z$ are lowest weight vectors of simple $(d+1)$-dimensional $\mathfrak s$-modules, so that
\[
E^{d/2}x=\frac1{d!}F^{d/2}E^dx,\qquad E^{d/2}z=\frac1{d!}F^{d/2}E^dz.
\]
Hence
\begin{multline*}
[E^{d/2}x,E^dz]=\frac1{d!}[F^{d/2}E^dx,E^dz]=(-1)^{d/2}\frac1{d!}[E^dx,F^{d/2}E^dz]\\
=(-1)^{d/2}[E^dx,E^{d/2}z].
\end{multline*}
\end{proof}

\begin{Proposition}
Any Cartan subspace for the representation of $Z_G(\sa)$ in $\g_{-d}$ is a subalgebra with respect to the product $*$. Hence it is called a Cartan subalgebra.
\end{Proposition}
\begin{proof}
Let $\qa$ be a minimal reducing subalgebra. Then $\qa_{-d}$ is a subalgebra of $(\g_{-d},*)$ and a Cartan subspace for $Z_G(\sa)|\g_{-d}$.
\end{proof}

\begin{Corollary}
All Cartan subalgebras in the algebra $(\g_{-d},*)$ are conjugate.
\end{Corollary}
\begin{proof}
It follows from Theorem \ref{cartanweyl} (a) and Proposition \ref{proporpolar}.
\end{proof}

Note that in the case when $d/2$ is even we get the usual Cartan subalgebras. In the case when $d/2$ is odd and $e$ is not irreducible, then either $\g_{-d}=\qa_{-d}$ for a minimal reducing subalgebra $\qa$, or we get Cartan subalgebras in simple Jordan algebras, which can be defined as maximal associative semisimple subalgebras. Their conjugacy is discussed in \cite{J}.

Now we turn to the identification of the algebras $(\g_{-d},*)$, defined by \eqref{operdef}, as listed in Tables \ref{irreds}, \ref{reductions}, \ref{reductionse6}, \ref{reductionse7}, \ref{reductionse8} and \ref{reductionsf4}. We use the following properties of these algebras, which are either obvious or proved above:
\begin{enumerate}[label=(\alph*)]
\item\label{prodinvar} The product $*$ is $Z(\sa)$-invariant.
\item\label{forminvar} The space $\g_{-d}$ carries a non-degenerate symmetric $Z(\sa)$-invariant bilinear form $(\cdot,\cdot)$, which is associative for the product $*$.
\item\label{prodcomm} The product $*$ is commutative if $d/2$ is odd, and anticommutative if $d/2$ is even.
\item The representation of $Z(\sa)$ is a direct sum of at most two irreducible representations, provided that $\g$ is simple.
\item\label{redsubmd} For any reducing subalgebra $\qa\subset\g$ the subspace $\qa\cap\g_{-d}$ is a subalgebra of the algebra $(\g_{-d},*)$.
\end{enumerate}

The following two lemmas are useful for the identification of the product $*$ when $d/2$ is odd, resp. even.

\begin{Lemma}\label{unitlemma}
Let $({\mathfrak a},*)$ be a finite-dimensional unital commutative algebra with a non-degenerate associative symmetric bilinear form $(\cdot,\cdot)$, invariant with respect to a group $G$ of automorphisms of $\mathfrak a$. Suppose that, with respect to the group $G$, $\mathfrak a$ decomposes as a trivial $1$-dimensional and non-trivial irreducible representation $V$, with $V*V\nsubseteq\field\bf1$, and that there is a unique, up to a scalar factor, map of $G$-modules ${\mathrm S}^2V\to V$. Then such a product on $\mathfrak a$ is unique, up to isomorphism.
\end{Lemma}
\begin{proof}
Note that ${\mathfrak a}=\field{\bf1}\oplus V$ is the decomposition of $\mathfrak a$ in an orthogonal direct sum of $G$-invariant subspaces and that the bilinear form can be normalized in such a way that $({\bf1},{\bf1})=1$. For $a,b\in\mathfrak a$ write $a*b=\overline{a*b}+\alpha\bf1$, where $\alpha\in\field$, $\overline{a*b}\in V$. Then, taking inner product with $\bf1$ and using associativity of the bilinear form, we obtain:
$$
({\bf1},a*b)=(a,b)=(\overline{a*b}+\alpha{\bf1},{\bf1})=\alpha.
$$
Hence $a*b=\overline{a*b}+(a,b)\bf1$.
\end{proof}

\begin{Lemma}\label{zerolemma}
Let $({\mathfrak a},[-,-])$ be a finite-dimensional skew-commutative algebra with a non-degenerate associative symmetric bilinear form $(\cdot,\cdot)$, invariant with respect to a group $G$ of automorphisms of $\mathfrak a$. Suppose that, with respect to the group $G$, $\mathfrak a$ decomposes as a trivial $1$-dimensional and non-trivial irreducible representation $U$, with $[U,U]\ne0$. Suppose that there exists a unique, up to a scalar factor, map of $G$-modules $\Lambda^2U\to U$. Then such a product on $\mathfrak a$ is unique, up to isomorphism.
\end{Lemma}
\begin{proof}
As in the previous lemma, we may assume that $({\bf1},{\bf1})=1$, $(U,{\bf1})=0$, and that restriction of $(\cdot,\cdot)$ to $U$ is nondegenerate.
For $a,b\in U$, write $[a,b]=p[a,b] +\alpha{\bf1}$, where $\alpha\in\field$ and $p$ is the projection on $U$.
Taking inner product with $\bf1$, we get $([a,b],{\bf1})=\alpha$, in particular, $\alpha=0$ if $a=b$.
We have $[b,{\bf1}]=\beta b$ for $b\in U$, with $\beta\in\field$ independent of $b$.
Then due to associativity of the form, $(a,[b,{\bf1}])=\alpha$, hence $\beta (a,b)=\alpha$.
Taking $a=b$ we obtain, as above, $\alpha=0$, hence $\beta(a,a)=0$.
Since $(\cdot,\cdot)$ is non-degenerate on $U$, we conclude that $\beta=0$.
Hence $[U,U]\subseteq U$, and $\g_{-d}$ is a direct sum of the algebra $U$ and a trivial 1-dimensional algebra ${\mathbb F}{\bf1}$.
Since on $U$ the product is non-zero and up to a scalar there is a unique $G$-invariant linear map $\Lambda^2U\to U$, we conclude that the product on $U$ is uniquely defined up to
a non-zero scalar.
\end{proof}

\

\noindent{\bf5B.}
Lemmas \ref{unitlemma} and \ref{zerolemma} are used in order to identify the algebra structure $(\g_{-d},*)$ in cases when $d/2$ is odd and even respectively. The lemmas are not applicable only in a few cases of nilpotent elements in exceptional Lie algebras, when the result can be checked directly on the computer. In many cases the algebras $(\g_{-d},*)$ are isomorphic to the well-known Lie or Jordan algebra structures; however in general they are neither Lie nor Jordan.

General nonassociative commutative algebras have been studied by various authors --- see e. g. \cite{walcher} (and many others). Much information about their appearance in connection with various questions of differential geometry has been provided in \cite{MOFox}.

All *-algebras that appear for irreducible nilpotent elements with odd $d/2$ fall into the series of algebras $\CA_\lambda(n)$ with the basis $p_1$, ..., $p_n$ that have multiplication table
$$
p_i^2=p_i,\quad p_ip_j=\lambda(p_i+p_j),\qquad 1\le i\ne j\le n, \lambda\in\field.
$$
For most $\lambda$, these algebras are not Jordan --- in fact, they are Jordan only for $n=1$, or $n\ge2$ and $\lambda=\frac12$, or $\lambda=0$ (in the latter case they are associative).

On the other hand, it is easy to check that all algebras $\CA_\lambda(n)$ satisfy two quartic identities. Namely, denoting by $\abr{x,y,z}=(xy)z-x(yz)$ the associator, every $a,b,c,d\in\CA_\lambda(n)$ satisfy
\begin{equation}\label{confass}
\abr{a,b,c}d-\abr{a,d,c}b=(ab)(cd)-(ad)(bc)
\end{equation}
and
\begin{equation}\label{mockjordan}
\abr{a,bd,c}+\abr{b,cd,a}+\abr{c,ad,b}=0.
\end{equation}
The identity \eqref{mockjordan} can be also equivalently written in terms of the multiplication operators $L_x$, i.~e. the operators given by $L_x(y)=xy$:
$$
[L_a,L_b]L_c+[L_b,L_c]L_a+[L_c,L_a]L_b=0.
$$
Note close resemblance to the Jordan identity, which is equivalent to
$$
\abr{ab,d,c}+\abr{bc,d,a}+\abr{ca,d,b}=0,
$$
or in terms of the multiplication operators,
$$
[L_{ab},L_c]+[L_{bc},L_a]+[L_{ca},L_b]=0.
$$

As pointed out by V.~Sokolov \cite{sokolov}, the identity \eqref{mockjordan} is actually a consequence of the identity \eqref{confass}.

For $\lambda\ne\frac12$ (which is the case for all of our irreducible nilpotent elements) the algebra $\CA_\lambda(n)$ has finitely many idempotents; since the equations determining idempotency are quadratic and there are $n$ of them, by B\'ezout's theorem the number of nonzero idempotents is less than $2^n$. In fact, $\frac1{2k\lambda+1}\sum_{i\in S}p_i$ is an idempotent of $\CA_\lambda(n)$ for any subset $S$ of $\{1,...,n\}$ of cardinality $k+1$. For $\lambda=-\frac1{2k}$ with integer $0<k<n-1$ this gives  $2^n-\binom n{k+1}$ idempotents, while for all other $\lambda\ne\frac12$, $\CA_\lambda(n)$ has exactly $2^n-1$ distinct nonzero idempotents. This is the case in all of our situations too, so that our *-algebras with $n$-dimensional $\g_{-d}$ have $2^n-1$ distinct 1-dimensional subalgebras.

It is clear from the multiplication table that the subspace of $\CA_\lambda(n)$ spanned by any subset $S\subseteq\{p_1,...,p_n\}$ is a subalgebra (isomorphic to $\CA_\lambda(k)$, where $k$ is the cardinality of $S$). Further subalgebras can be obtained from these via actions by algebra automorphisms. While there is an obvious action of S$_n$ through permuting the generators $p_i$, there are no other apparent automorphisms except for $\lambda=-\frac1{n-1}$: indeed, in this case $p_0=-p_1-...-p_n$ is an idempotent and moreover  $p_0p_i=\lambda(p_0+p_i)$, so that there will be additional automorphisms permuting $p_0$ with all other $p_i$. Thus the automorphism group of $\CA_{-\frac1{n-1}}(n)$ contains S$_{n+1}$. As shown in \cite{harada}, $\CA_{-\frac1{n-1}}(n)$ does not have any further automorphisms, so that its automorphism group is exactly S$_{n+1}$ (cf. the last two lines of Table \ref{irreds}).

In the cases occurring in Table \ref{reductions} we can explicitly describe all subalgebras of $\CA_\lambda(n)$. For $\CA_\lambda(2)$ every proper subalgebra is 1-dimensional. For $\CA_\lambda(3)$, looking directly at the conditions on a 2-dimensional subspace to be a subalgebra, we obtain that for $\lambda\ne\frac12$ there are exactly six 2-dimensional subalgebras, namely, those spanned by $\abr{p_1,p_2}$, $\abr{p_1,p_3}$, $\abr{p_2,p_3}$, $\abr{p_1+p_2,p_3}$, $\abr{p_1+p_3,p_2}$ and $\abr{p_2+p_3,p_1}$. Similarly, for $\lambda\ne\frac12$ the algebra $\CA_\lambda(4)$ has only ten 3-dimensional subalgebras $\abr{p_i,p_j,p_k}$ and $\abr{p_i,p_j,p_k+p_\ell}$; and for the only 4-dimensional case $\CA_{-\frac13}(4)$ in Table \ref{reductions} there are only 25 2-dimensional subalgebras $\abr{p_i,p_j}$, $\abr{p_i+p_j,p_k}$, $\abr{p_i+p_j+p_k,p_\ell}$ and $\abr{p_i+p_j,p_k+p_\ell}$, for pairwise distinct $i$, $j$, $k$, $\ell$. It thus follows that for algebras $\CA_\lambda(n)$ occurring in Table \ref{reductions}, all subalgebras are spanned by idempotents.

First consider the case when $e$ is an irreducible nilpotent element. It follows from Table \ref{irreds} that for $d/2$ even we always have $\dim\g_{-d}=1$, and since the product $*$ is anticommutative, the algebra $(\g_{-d},*)$ has zero multiplication. Next, when $d/2$ is odd and $\g$ is an exceptional Lie algebra, we identify the algebra $(\g_{-d},*)$ with the aid of computer, as follows.

Structure constants table of $(\g_{-d},*)$ in the root vector basis, and the $\z(\sa)$-module structure are computed using GAP. When the algebra is commutative, in each case idempotents are computed using a generic element of $\g_{-d}$ with indeterminate coefficients. When there is a basis consisting of idempotents, the algebra is identified with one of the $\CA_\lambda(n)$ using it. In Appendix \ref{appB}, the cases $J_n^c(BD)$ are identified finding a basis with almost all pairwise products zero, and the algebras $H_5$ and $H_8$ are identified using explicit isomorphisms.

The cases when $(\g_{-d},*)$ is a Lie algebra are determined using the command $\mathtt{TestJacobi}$, and then the isomorphism type of this algebra is determined using the commands $\mathtt{LeviMalcevDecomposition}$ and $\mathtt{SemiSimpleType}$ in \cite{gap}.

Finally there are cases when the product is skew-commutative and does not satisfy the Jacobi identity. These cases are identified with the 7-dimensional simple Malcev algebra.

The remaining cases are treated by the following two lemmas.

\begin{Lemma}
In cases $2_k$ and $3_k$ of Table \ref{irreds} the product $*$ is non-zero.
\end{Lemma}
\begin{proof}
For $e$ the principal nilpotent element and $F$ the lowest root vector --- in $\spa(2k)$ for the case 2$_k$ and in $\so(2k+1)$ for the case 3$_k$ --- according to \eqref{operdef} we have to show that the element
$$
[(\ad e)^{2k-1}F,F]
$$
is nonzero.

Recall the well-known identity in any associative algebra (see e.~g. \cite{IDLA}*{(3.8.1)}):
$$
\exp(a)b\exp(-a)=(\exp(\ad a))b.
$$
Using this identity in the standard representation we have that $(\ad e)^jF$ is a scalar multiple of the coefficient at $t^j$ of the matrix $\exp(te)F\exp(-te)$. More precisely,
\begin{equation}\label{adeformula}
(\ad e)^jF=\sum_{i=0}^j(-1)^i\binom jie^{j-i}Fe^i.
\end{equation}

For the case $2_k$, in the standard representation on $\field^{2k}$ the matrix for $e$ is the largest Jordan block, while the only nonzero entry of the matrix for $F$ is 1 in the lower left corner. It follows that the coefficient at $t^{2k-1}$ of $\exp(te)F\exp(-te)$ is the diagonal matrix with entries $\frac{(-1)^{i-1}}{(i-1)!(2k-i)!}$, $i=1,...,2k$. Moreover, for a diagonal matrix $D$, the matrix $[D,F]$ is $-D_{1,1}+D_{2k,2k}$ times $F$. In our case these diagonal entries have equal absolute values and opposite signs, so that this gives $-\frac2{(2k-1)!}F\ne0$.

For the case $3_k$, in the standard representation on $\field^{2k+1}$ there also is a basis such that the matrix of $e$ is the largest Jordan block. In this basis, the matrix for the lowest root vector $F$ has $(-1)^k2$ at positions $(2k,1)$ and $(2k+1,2)$ and zeroes elsewhere. Thus for any diagonal matrix $D$ the matrix $[D,F]$ has $(-1)^k2(-D_{1,1}+D_{2k,2k})$ at the $(2k,1)$st position, $(-1)^k2(-D_{2,2}+D_{2k+1,2k+1})$ at the $(2k+1,2)$nd position and zeroes elsewhere.

Moreover the coefficient at $t^{2k-1}$ of $\exp(te)F\exp(-te)$ is the diagonal matrix $D$ with entries
$$
D_{i,i}=
(-1)^k2\begin{cases}
\frac1{(2k-1)!},&i=1,\\
\frac{(-1)^{i-1}}{(i-1)!(2k-i)!}+\frac{(-1)^{i-2}}{(i-2)!(2k+1-i)!},&1<i<2k+1,\\
-\frac1{(2k-1)!},&i=2k+1.
\end{cases}
$$
It follows that $[D,F]$ in this case is $(-1)^k\frac{2(2k-3)}{(2k-1)!}F\ne0$.
\end{proof}

\begin{Lemma}\label{case4}
In cases $4_k$ the algebra $(\g_{-d},*)$ is isomorphic to the algebra $\CA_{-k}(2)$.
\end{Lemma}
\begin{proof}
We will use the same basis of the standard representation that was described in \eqref{basis4k}, with the choice of $e$ such that it acts on this basis as indicated there:
$$
\setlength\arraycolsep{2pt}
\begin{array}{lllllllllllllllllll}
x_{-k-1}&\mapsto&x_{-k}&\mapsto&\cdots&\mapsto&x_{-1}&\mapsto&x_0&\mapsto&x_1&\mapsto&\cdots&\mapsto&x_k&\mapsto&x_{k+1}&\mapsto&0,\\
      &       &y_{-k}&\mapsto&\cdots&\mapsto&y_{-1}&\mapsto&y_0&\mapsto&y_1&\mapsto&\cdots&\mapsto&y_k&\mapsto&0,
\end{array}
$$
i.~e. the matrix of $e$ in the standard representation consists of two Jordan blocks, of sizes $2k+3$ and $2k+1$. It will be convenient for us to choose $F_1,F_2$ in such a way that $(F_1,(-1)^kF_2)$ is the root vector basis of $\g_{-d}$, with $F_1$ the lowest root vector. In the above basis of the standard representation these then act as follows:
$$
F_1(x_{k+1})=y_{-k}, F_1(y_k)=-x_{-k-1}, F_2(x_{k+1})=x_{-k}, F_2(x_k)=x_{-k-1},
$$
both sending all remaining basis elements to zero.

We will compute the multiplication table of $\g_{-d}$ in this basis, i.~e. find
$$
F_i*F_j=[(\ad e)^{2k+1}F_i,F_j],\quad i,j=1,2.
$$
Let then $\bar F_1=(\ad e)^{2k+1}F_1$, $\bar F_2=(\ad e)^{2k+1}F_2$. Using again \eqref{adeformula}, we find
\begin{align*}
&\bar F_1(x_{k+1})=\bar F_1(x_{-k-1})=0,\\
&\bar F_1(x_j)=-(-1)^{k+j}\binom{2k+1}{k+j}y_j, \bar F_1(y_j)=-(-1)^{k-j}\binom{2k+1}{k-j}x_j, -k\le j\le k
\intertext{and}
&\bar F_2(x_{k+1})=x_{k+1}, \bar F_2(x_{-k-1})=-x_{-k-1},\\
&\bar F_2(x_j)=-(-1)^{k-j}\left(\binom{2k+1}{k+j}-\binom{2k+1}{k-j}\right)x_j, \bar F_2(y_j)=0, -k\le j\le k.
\end{align*}
From this we get the multiplication table,
$$
F_1*F_1=-(2k+1)F_2,\qquad F_1*F_2=F_2*F_1=-F_1,\qquad F_2*F_2=(2k-1)F_2.
$$

By solving $(\alpha_1F_1+\alpha_2F_2)*(\alpha_1F_1+\alpha_2F_2)=\alpha_1F_1+\alpha_2F_2$ for $\alpha_1$, $\alpha_2$, we find that the elements
$$
P=-\frac{F_1+F_2}2, Q=\frac{F_1-F_2}2
$$
are idempotents, and moreover
\begin{multline*}
P*Q=\left(-\frac{F_1+F_2}2\right)*\left(\frac{F_1-F_2}2\right)=\frac{F_2*F_2-F_1*F_1}4\\
=\frac{(2k-1)F_2+(2k+1)F_2}4=kF_2=-k(P+Q),
\end{multline*}
which gives the multiplication table for $\CA_{-k}(2)$.
\end{proof}

In all irreducible cases one has

\begin{Proposition}
If $e$ is irreducible, then for any $F\in\g_{-d}$, the cyclic element $e+F$ is semisimple if and only if $F$ does not lie in any proper subalgebra of $(\g_{-d},*)$.
\end{Proposition}
\begin{proof}
This is clear when $\dim\g_{-d}=1$. For the case 4$_k$ this follows by comparing computations with \eqref{basis4k} and the proof of Lemma \ref{case4} above.
Indeed with the former we saw that, for some particular choice of $e$, the element $e+\lambda_1F_1+\lambda_2F_2$ is semisimple if and only if $\lambda_1\ne0$ and $\lambda_2\ne\pm\lambda_1$, where $(F_1,F_2)$ is the root vector basis of $\g_{-d}$, with $F_1$ the lowest root vector. While with the latter, for the same choice of $e$, we saw that nonzero idempotents in the algebra $(\g_{-d},*)$ are $-\frac{F_1+F_2}2$, $\frac{F_1-F_2}2$ and $\frac{F_2}{2k-1}$, so that there are three proper subalgebras, spanned by these elements. But these are precisely 1-dimensional subspaces spanned by an element $\lambda_1F_1+\lambda_2F_2$ with $\lambda_1=-\lambda_2$, $\lambda_1=\lambda_2$ and $\lambda_1=0$ respectively.

In the remaining cases of irreducible $e$ (cases 7,11,16,17,18 of Table \ref{irreds}) we similarly compare semisimplicity condition on a generic cyclic element with the algebra structure on $(\g_{-d},*)$. As an illustration, let us treat here the last of these cases, 18 (nilpotent element with label E$_8(a_7)$, depth 10, $\dim(\g_{-10})=4$) --- other cases are similar but shorter. Let us choose an orbit representative $e$ in the form
$$
e_{\foreight{.5}00001000}+e_{\foreight{.5}01111000}+e_{\foreight{.5}01011110}+e_{\foreight{.5}11121000}+
e_{\foreight{.5}11111100}+e_{\foreight{.5}10111110}+e_{\foreight{.5}01121100}+e_{\foreight{.5}00111111}.
$$
Let
$$
F_1=f_{\foreight{.5}23465321},\ F_2=f_{\foreight{.5}23465421},\ F_3=f_{\foreight{.5}23465431},\ F_4=f_{\foreight{.5}23465432}
$$
be the root vector basis of $\g_{-10}$. As explained before (at the start of the proof of \eqref{conj}), it follows from \cite{springer}*{9.5} that a cyclic element $C=e+\lambda_1F_1+\lambda_2F_2+\lambda_3F_3+\lambda_4F_4$ is semisimple if and only if it is regular semisimple. Then regular semisimplicity can be checked by looking at the appropriate coefficient of the characteristic polynomial for $\ad C$. In our case this coefficient turns out to be a scalar multiple of a power of
\begin{equation}\label{charp}
\begin{aligned}
&(\lambda_1+\lambda_2-\lambda_3-\lambda_4) (\lambda_1-\lambda_2+\lambda_3-\lambda_4) (-\lambda_1+\lambda_2+\lambda_3-\lambda_4)\\
&\qquad\times (\lambda_1+\lambda_2+\lambda_3+\lambda_4) (\lambda_1^2-\lambda_2^2)(\lambda_1^2-\lambda_3^2)(\lambda_2^2-\lambda_3^2).
\end{aligned}
\end{equation}
On the other hand, computing $F_i*F_j=[(\ad e)^5F_i,F_j]$ gives
\begin{multline*}
F_1*F_2=10F_3,\, F_2*F_3=10F_1,\, F_1*F_3=10F_2,\, F_1*F_1=F_2*F_2=F_3*F_3=10F_4,\\
F_4*F_4=-6F_4, \text{ and } F_i*F_4=2F_i, i=1,2,3.
\end{multline*}
One checks that with respect to this multiplication the elements
\begin{multline*}
P_1=\frac{F_1-F_2-F_3+F_4}{24}, P_2=\frac{-F_1+F_2-F_3+F_4}{24}, P_3=\frac{-F_1-F_2+F_3+F_4}{24},\\
P_4=\frac{F_1+F_2+F_3+F_4}{24}
\end{multline*}
are idempotents and satisfy
$$
P_i*P_j=-\frac13(P_i+P_j), i,j\in\{1,2,3,4\}, i\ne j.
$$
It follows that $(\g_{-10},*)$ is isomorphic to $\CA_{-\frac13}(4)$ and its maximal (3-dimensional) subalgebras are spanned by linearly independent triples from the set of vectors $P_i$, $P_i+P_j$, $P_i+P_j+P_k$, $P_1+P_2+P_3+P_4$. This amounts to ten 3-dimensional subspaces, four spanned by $\{P_i,P_j,P_k\}$, $\{i,j,k\}\subset\{1,2,3,4\}$ and six spanned by $\{P_i,P_j,P_1+P_2+P_3+P_4\}$, $\{i,j\}\subset\{1,2,3,4\}$. It is then straightforward to check that the subspace spanned by $\{P_1,P_2,P_3\}$ consists of $\lambda_1F_1+\lambda_2F_2+\lambda_3F_3+\lambda_4F_4$ with $\lambda_1+\lambda_2+\lambda_3+\lambda_4=0$, that spanned by $\{P_i,P_j,P_4\}$, $\{i,j\}\subset\{1,2,3\}$ corresponds to $\lambda_i+\lambda_j=\lambda_k+\lambda_4$, with $\{k\}=\{1,2,3\}\setminus\{i,j\}$, the one spanned by $\{P_i,P_j,P_1+P_2+P_3+P_4\}$, $\{i,j\}\subset\{1,2,3\}$ corresponds to $\lambda_i+\lambda_j=0$, and the one spanned by $\{P_i,P_4,P_1+P_2+P_3+P_4\}$, $i\in\{1,2,3\}$ corresponds to $\lambda_j=\lambda_k$, where $\{j,k\}=\{1,2,3\}\setminus\{i\}$. Comparing these to \eqref{charp} we see that indeed $C=e+\lambda_1F_1+\lambda_2F_2+\lambda_3F_3+\lambda_4F_4$ loses semisimplicity if and only if $\lambda_1F_1+\lambda_2F_2+\lambda_3F_3+\lambda_4F_4$ belongs to a proper subalgebra of $(\g_{-10},*)$.
\end{proof}

\

\noindent{\bf5C.}
We return to the identification of the algebra $(\g_{-d},*)$ with those listed in Tables \ref{reductions}, \ref{reductionse6}, \ref{reductionse7}, \ref{reductionse8} and \ref{reductionsf4}.

Recall that a Malcev algebra is defined by a skewsymmetric bracket, satisfying a quartic identity, which is implied by the Jacobi identity (thus any Lie algebra is a Malcev algebra). It was proved in \cite{sagle} and \cite{kuzmin} that any simple finite-dimensional Malcev algebra is either one of the simple Lie algebras, or is the 7-dimensional space of imaginary octonions, equipped with the usual bracket $[a,b]=ab-ba$. We denote the latter algebra by $M$.

Recall that isomorphism classes of simple finite-dimensional Jordan algebras are in bijective correspondence to conjugacy classes of even nilpotent elements $e$ of depth $2$ in simple Lie algebras (see  \cite{J}). Namely the product $*$  on $\g_{-2}$ defines a structure of a Jordan algebra and all simple Jordan algebras are thus obtained. The complete list consists of all $n\times n$ matrices with product $a\cdot b=ab+ba$, which we denote by $J_n(A)$, the subalgebra of $J_n(A)$ consisting of matrices selfadjoint with respect to a symmetric (respectively skewsymmetric) non-degenerate bilinear form, which we denote by $J_n(C)$ (resp. $J_n(D)$), and the space $V\oplus\field1$, where $V$ is the $n$-dimensional space with a non-degenerate symmetric bilinear form $(\cdot\mid\cdot)$, with product $a\cdot b=(a\mid b)1$, $a\cdot1=1\cdot a=a$ for $a,b\in V$, $1\cdot1=1$, which we denote by $J_n(BD)$. Finally there is the 27-dimensional exceptional Albert's algebra which we denote by $J(E)$. All of these Jordan algebras are simple. This notation stems from the fact that these Jordan algebras correspond to nilpotent elements in the Lie algebras of the corresponding type A, B, C, D or E$_7$.

It suffices to identify the algebra $(\g_{-d},*)$ in the cases $\g=\g^{\mathrm{ev}}$, using the passage from $\g$ to $\g^{\mathrm{ev}}$, described in Appendix \ref{appA}. The ``shortest'' case $d=2$ of $e\in\g^{\mathrm{ev}}$ corresponds to an even nilpotent element of depth 2. As mentioned above, conjugacy classes of these nilpotent elements correspond bijectively to the isomorphism class of a structure of a simple Jordan algebra on $\g_{-2}$.

Next, consider the case $d/2$ odd and $>1$. By property \ref{redsubmd}, for the nilpotent elements $e$ with $\qa_{\max}=\qa_{\min}$ the identification of $(\g_{-d},*)$ reduces to that of $(\qa_{\min},*)$, which is the case of irreducible nilpotent elements, discussed above. As a result, only the following nilpotent elements with $d/2$ odd remain to be considered:
$$
\sla(2kn)\ni[(2k)^{(n)}],\ \spa(2kn)\ni[(2k)^{(n)}],\ \so(4kn)\ni[(2k)^{(2n)}].
$$
But in all these cases $\qa_{\min}$ is a commutative associative semisimple subalgebra and the representation of $Z(\sa)$ on $\g_{-d}$ is a direct sum of a non-trivial irreducible and the trivial 1-dimensional subrepresentations. This and properties \ref{prodinvar}, \ref{forminvar}, \ref{prodcomm} along with Lemma \ref{unitlemma} allow us to identify the algebras $(\g_{-d},*)$ with the Jordan algebras $J_n(A)$, $J_n(C)$ and $J_n(D)$ respectively. In order for Lemma \ref{unitlemma} to be applicable here requires ensuring that the symmetric squares of the representations ad$_{\sla(n)}$, S$^2($st$_{\so(n)})$ and $\Lambda^2(\text{st}_{\spa(2n)})$ each contain a unique copy of the same representation, respectively. It can be checked e.~g. using \cite{OV}*{Table 5 (pages 300--303)}. The least obvious of these cases is the one for $J_n(D)$. In this case $n=2j$ is even; consider the involution $\iota$ on the algebra of $2j\times2j$ matrices given by
\begin{equation}\label{invol}
\iota\left(\begin{array}{rr}A_{11}^{\phantom\intercal}&A_{12}^{\phantom\intercal}\\A_{21}^{\phantom\intercal}&A_{22}^{\phantom\intercal}\end{array}\right):=\left(\begin{array}{rr}A_{22}^\intercal&-A_{12}^\intercal\\-A_{21}^\intercal&A_{11}^\intercal\end{array}\right).
\end{equation}
Fixed points of this involution consist of $2\times2$ blocks of $j\times j$ matrices with skew-symmetric $A_{12}$ and $A_{21}$ and with $A_{22}=A_{11}^\intercal$. They thus can be identified with the exterior square of a $2j$-dimensional space through the canonical isomorphism $\Lambda^2(V^*\oplus V)\cong\Lambda^2(V)^*\oplus\gl(V)\oplus\Lambda^2(V)$ for a $j$-dimensional space $V$. They are closed under anticommutator and form a simple Jordan algebra of symplectic type, acted upon via derivations by commutators with the Lie algebra of anti-fixed points of $\iota$. The latter in turn can be identified with the symmetric square of a $2j$-dimensional space through S$^2(V^*\oplus V)\cong$~S$^2(V)^*\oplus\gl(V)\oplus$~S$^2(V)$, being blocks with $A_{22}=-A_{11}^\intercal$ and $A_{12}$, $A_{21}$ symmetric, which is the Lie algebra $\spa(2j)$, with respect to the standard skew-symmetric form $\omega$ on $V^*\oplus V$ given by $\omega(\ph,\ph')=\omega(v,v')=0$, $\omega(\ph,v)=-\omega(v,\ph)=\ph(v)$.

It remains to consider the case when $d/2$ is even. As before, when $\dim\g_{-d}=1$ we have the 1-dimensional algebra with zero multiplication. Since we may assume that $\g=\qa_{\max}$, we are left with the following cases:
$$
\sla(2kn+n)\ni[(2k+1)^{(n)}];\ \spa(4kn+2n)\ni[(2k+1)^{(2n)}];\ \so(4kn+2n)\ni[(2k+1)^{(2n)}]
$$
for classical Lie algebras, and the following cases for exceptional Lie algebras:
$$
{\mathrm E}_6\ni[2{\mathrm A}_2];\ {\mathrm E}_7\ni[2{\mathrm A}_2], [{\mathrm A}_6];\ {\mathrm E}_8\ni[2{\mathrm A}_2], [{\mathrm A}_6];\ {\mathrm F}_4\ni[\tilde{\mathrm A}_2].
$$
In all these cases there exists a unique, up to constant factors, product, satisfying properties \ref{prodinvar}, \ref{forminvar}, \ref{prodcomm}. It remains to prove that product $*$ in these cases is non-zero on each non-trivial irreducible component of the $Z(\sa)$-module $\g_{-d}$. For $\g$ of exceptional type, this is done by direct calculation: with the GAP command $\mathtt{DirectSumDecomposition}$ the irreducible components are found, and the $*$-products of generic elements of these components are computed to be nonzero (as mentioned, we use the \texttt{SLA} package by W. de Graaf \cite{deGraaf} for the \texttt{GAP} system \cite{gap}). As an example, take the case $2{\mathrm A}_2$ in ${\mathrm E}_8$. Here $Z(\sa)$ is the semidirect product of a 2-element cyclic group with $\rG_2\times\rG_2$. Here $\g_{-d}$ is 14-dimensional and representation of $Z(\sa)$ on it realizes two copies of the 7-dimensional irreducible representation of $\rG_2$. Decomposing the exterior square of this representation we find that it contains two 7-dimensional irreducible components. Since the $*$-product must be $\rG_2\times\rG_2$-invariant, we deduce that each of these components can only map nontrivially to separate 7-dimensional summands in $\g_{-d}$. We then check by direct calculation that there are indeed nonzero products on each of these separately. We then finally conclude that the algebra structure is isomorphic to that of two copies of the simple 7-dimensional Maltsev algebra.

For $\g$ in Table \ref{reductions}, applicability of Lemmas \ref{unitlemma} and \ref{zerolemma} when $d/2$ is odd, resp. even, still requires to show that there is at least one instance of the $*$-product with nonzero projection to the nontrivial irreducible summand. This follows from

\begin{Lemma}
Let $e$ be a nilpotent element with partition of the form $[n^{(m)}]$ in a classical simple Lie algebra $\g$. Then the algebra $(\g_{-d},*)$ is as in Table \ref{reductions}.
\end{Lemma}
\begin{proof}
We can choose a basis in the standard representation in such a way that elements of $\g$ are represented by block matrices, consisting of $n\times n$ blocks of size $m\times m$ each, in such a way that in this basis $e$ is ``block-principal'', i.~e. represented by a matrix with identity matrices in blocks $E_{12}$, $E_{23}$, ..., $E_{n-1,n}$ and zeroes elsewhere, while elements $F\in\g_{-d}$ are represented by a single block $F_{n1}=F$ with zeroes elsewhere. Moreover, using the argument from \cite{ekv}*{Section 4}, this basis can be chosen in such a way that the $m\times m$ matrix $F_{n1}$ is
\begin{itemize}
\item symmetric if $\g=\spa(mn)$ with $n$ even,
\item anti-fixed point of the involution \eqref{invol} if $\g=\spa(mn)$ with $n$ odd (hence $m$ even),
\item fixed point of the involution \eqref{invol} if $\g=\so(mn)$ with $n$ even (hence $m$ even),
\item skew-symmetric if $\g=\so(mn)$ with $n$ odd.
\end{itemize}

In our case $d=2n-2$, and  using \eqref{adeformula} we see that the matrix $(\ad e)^{\frac d2}F$ is block-diagonal, with matrices $D_{ii}=(-1)^{i-1}\binom{n-1}{i-1}F$, $i=1,...,n$, along the diagonal. Consequently
\begin{multline*}
F*F'=[(\ad e)^{n-1}F,F']=[\operatorname{diag}(D_{11},...,D_{nn}),F'_{n1}]=D_{11}F'-F'D_{nn}\\=FF'-(-1)^{n-1}F'F,
\end{multline*}
so that the algebra structure on $\g_{-d}$ is indeed as claimed.
\end{proof}

\ 

Examining the respective instances in Tables \ref{reductions}, \ref{reductionsf4}, \ref{reductionse6}, \ref{reductionse7}, \ref{reductionse8} we arrive at

\ 

\begin{Theorem}\label{algm}
There are the following three possibilities for a nilpotent element $e$ of semisimple type.
\begin{itemize}
\item[(a)] $\rank e=\dim\g_{-d}$ and $d/2$ is odd (resp. even). Then the algebra $\g_{-d}$ with product \eqref{operdef} is isomorphic to one of the commutative algebras $\CA_\lambda(n)$, where $n=\dim\g_{-d}$ (resp. to the $1$-dimensional Lie algebra);
\item[(b)] $\rank e<\dim\g_{-d}$ and $d/2$ is odd. Then the algebra $\g_{-d}$ with product \eqref{operdef} is isomorphic to one of the simple Jordan algebras;
\item[(c)] $\rank e<\dim\g_{-d}$ and $d/2$ is even. Then the algebra $\g_{-d}$ with product \eqref{operdef} is isomorphic to a direct sum of at most two simple Malcev algebras, including the $1$-dimensional one.
\end{itemize}
\end{Theorem}

As explained in the introduction, by looking at the tables, we obtain the following theorem.
\begin{Theorem}\label{theoS} Let $e$ be a nilpotent element of semisimple type in a simple Lie algebra~$\g$. We have the following
description of the set $\Ss_\g(e)$:
\begin{itemize}
\item[]Case (a) of Theorem \ref{algm}: $F$ lies outside of the union of hyperplanes, spanned by idempotents.
\item[]Case (b) of Theorem \ref{algm}:
\begin{itemize}
\item[(i)] $(\g_{-d},*)\not\simeq J_n(BD)$, then
\[
\ \hskip-10.5ex\Ss_\g(e)=\setof{F\in\g_{-d}}{\text{$Z_G(\sa)F$ is closed and $L_F\in\operatorname{End}\g_{-d}$ has maximal rank}},
\]
\item[(ii)] $(\g_{-d},*)\simeq J_n(BD)$, then
\[
\Ss_\g(e)=\setof{F\in\g_{-d}}{F*F\notin\field F}.
\]
\end{itemize}
\item[]Case (c) of Theorem \ref{algm}:
\[
\ \hskip-2.5ex\Ss_\g(e)=\setof{F\in\g_{-d}}{\text{$Z_\g(\sa)F$ is closed and $L_F\in\operatorname{End}\g_{-d}$ has maximal rank}}.
\]
\end{itemize}
\end{Theorem}

\begin{Conjecture}
Description of $F\in\Ss_\g(e)$, for which $G(e+F)$ has maximal dimension:
\begin{itemize}
\item[(i)] in cases (a) and b(ii) of Theorem \ref{theoS} all have maximal dimension,
\item[(ii)] in the remaining cases the orbit of $F$ has maximal dimension among the $Z_G(\sa)$-orbits in $\g_{-d}$, and $F\in\ca$, a Cartan subalgebra of $(\g_{-d},*)$,  lies outside of the union of reflection hyperplanes of the Weyl group of the polar linear group $Z_G^\circ(\sa)|\g_{-d}$.
\end{itemize}
\end{Conjecture}

\begin{Remark}
Let $e\in\g$ be a nilpotent element of even depth $d$, not divisible by $4$, and assume that $\dim\g_{-d}=1$, so that $\g_{-d}=\field a$ for some non-zero element $a$.
Then by \eqref{operdef} we have
$$
a*a=[(\ad e)^{\frac d2}(a),a]=P(e)a,
$$
where $P(e)$ is a homogeneous polynomial in $e\in\g_2$ of degree $\frac d2$. It is easy to see that this polynomial is $G_0$-semi-invariant, with character $\chi^{-1}$ where $\chi$ is the character for the action of $G_0$ on $\g_{-d}$. An interesting problem is to compute this polynomial. We found the answer in the case of a principal nilpotent element of a simple Lie algebra $\g$ of rank $r$. In this case $d=2(h-1)$, where $h$ is the Coxeter number. So $\frac d2$ is odd iff $h$ is even, which excludes $\g$ of type A$_n$, $n$ even. Write $e=\sum_{i=1}^rx_ie_i\in\g_2$, where $e_i$ are the root vectors attached to simple roots $\alpha_i$, and let $\theta=\sum_{i=1}^ra_i\alpha_i$ be the highest root. Then
$$
P(e)=\prod_{i=1}^rx_i^{a_i}.
$$
\end{Remark}

\appendix
\section*{Appendices}

\renewcommand{\thesubsection}{\Alph{subsection}}
\renewcommand{\thesubsubsection}{\Alph{subsection}.\arabic{subsubsection}}

\renewcommand{\theequation}{\thesubsection.\arabic{equation}}

\subsection{Even reductions}\label{appA}

Given a nilpotent element $e$ in a simple Lie algebra $\g$ with the standard $\sla(2)$-triple $\sa=\{e,f,h\}$, the even part $\g^{\mathrm{ev}}:=\bigoplus_k\g_{2k}$ of the grading of $\g$ is a subalgebra containing $e$, whose derived subalgebra is reducing, unless $e$ happens to be of nilpotent type (since then depth of $e$ in $\g^{\mathrm{ev}}$ drops by $1$). Denoting by $G$ the adjoint group of $\g$ and by $S\subset G^{\mathrm{ev}}$ the subgroups corresponding to $\sa$, resp.~$\g^{\mathrm{ev}}$, we easily see that $\g^{\mathrm{ev}}$ is the algebra of fixed points for an involution corresponding to the adjoint action of an order 2 element of $S$ which lies in the center of $G^{\mathrm{ev}}$.

Fixed point algebra of an order two inner automorphism of a simple Lie algebra $\g$ of rank $r$ is obtained by considering its extended Dynkin diagram whose nodes are labeled by coefficients $a_0=1$, $a_1$, ..., $a_r$ of the integer linear dependence of the columns of the extended Cartan matrix. Then a fixed point subalgebra of an inner involution is obtained by removing one node with label $2$ or two nodes with label $1$; in the second case one adds $T^1$ \cite{IDLA}*{Chapter 8}:

\

\begin{tabular}{ll}
algebra type&even subalgebra types\\
\hline
A$_n$&A$_k+$ A$_{n-k-1}+T^1$, $1\le k\le n-1$\\
B$_n$&B$_{n-1}+T^1$, D$_k+$ B$_{n-k}$, $2\le k\le n$\\
C$_n$&A$_{n-1}+T^1$, C$_k+$ C$_{n-k}$, $1\le k\le n-1$\\
D$_n$&A$_{n-1}+T^1$, D$_{n-1}+T^1$, D$_k+$ D$_{n-k}$, $2\le k\le n-2$\\
E$_6$&A$_5+$ A$_1$, D$_5+T^1$\\
E$_7$&D$_6+$ A$_1$, A$_7$, E$_6+T^1$\\
E$_8$&E$_7+$ A$_1$, D$_8$\\
F$_4$&B$_4$, C$_3+$ A$_1$\\
G$_2$&2A$_1$
\end{tabular}

\

\

For classical types, if $e$ is not even then the corresponding partition contains parts of both even and odd parities.
Let us separate this partition into two partitions, one containing even parts only and another odd parts only.
The derived subalgebra of $\g^{\mathrm{ev}}$ is the direct sum of two subalgebras, with $e$ decomposing into the sum of two nilpotent elements, one in each of these subalgebras,
with these two partitions. Here we assume that the partition with all parts equal to 1 corresponds to the zero nilpotent, i.~e. if the odd subpartition is such then $e$ has zero projection to the corresponding summand of $\g^{\mathrm{ev}}$.

\subsubsection{\bf Examples}
Let $e$ be a nilpotent element in $\g$ of type B$_8$ with partition $[5,2^{(4)},1^{(4)}]$. The odd subpartition $[5,1^{(4)}]$ is the partition of a nilpotent element in B$_4$ and the even one $[2^{(4)}]$ is the partition of a nilpotent element in D$_4$. Accordingly, $\g^{\mathrm{ev}}$ has type $\rB_4+\rD_4$, and $e\in\g^{\mathrm{ev}}$ decomposes into the sum of nilpotent elements with indicated partitions in these summands.

If $\g$ is of type C$_9$ and $e$ has partition $[4,2^{(2)},1^{(10)}]$, then the even subpartition $[4,2^{(2)}]$ belongs to a nilpotent element in C$_4$ and $[1^{(10)}]$ represents the zero nilpotent element in C$_5$. In this case $\g^{\mathrm{ev}}$ is C$_4+$ C$_5$, and $e$ belongs to the summand C$_4$, having partition $[4,2^{(2)}]$ there and projecting to zero in C$_5$.

\

For odd nilpotent elements in exceptional simple Lie algebras, we get the following picture. Nilpotent elements of nilpotent type are marked with an ``*''.

\setlength\tabcolsep{2pt}
\begin{tabular}[t]{rll}
\multicolumn3{l}{G$_2$, even subalgebra 2A$_1$:}\\
\\
&nilpotent $e$&partitions\\
\hline
&A$_1$\hfill\GDynkin{0,1}&$[2]$, $[1^{(2)}]$\bigstrut[t]\\
*&$\tilde{\mathrm A}_1$\hfill\GDynkin{1,0}&$[2]$, $[1^{(2)}]$
\end{tabular}%

\

\

\

\

\resizebox{.95\textwidth}{!}{%
\begin{tabular}{l|l}
\begin{tabular}[t]{rll}
\multicolumn3{l}{F$_4$, even subalgebra C$_3+$ A$_1$:}\\
\\
&nilpotent $e$&partitions in C$_3$, A$_1$\\
\hline
&A$_1$\hfill\FDynkin{0,1,0,0}&$[1^{(6)}]$, $[2]$\bigstrut[t]\\
*&A$_1+\tilde{\mathrm A}_1$\hfill\FDynkin{0,0,0,1}&$[2^{(3)}]$, $[1^{(2)}]$\\
*&$\tilde{\mathrm A}_2+$ A$_1$\hfill\FDynkin{1,0,0,1}&$[3^{(2)}]$, $[2]$\\
&C$_3(a_1)$\hfill\FDynkin{0,1,1,0}&$[4,2]$, $[1^{(2)}]$\\
&C$_3$\hfill\FDynkin{2,1,1,0}&$[6]$, $[1^{(2)}]$
\end{tabular}&
\begin{tabular}[t]{ll}
\multicolumn2{l}{F$_4$, even subalgebra B$_4$:}\\
\\
nilpotent $e$&partition in B$_4$\\
\hline
$\tilde{\mathrm A}_1$\hfill\FDynkin{1,0,0,0}&$[3,1^{(6)}]$\bigstrut[t]\\
A$_2+\tilde{\mathrm A}_1$\hfill\FDynkin{0,0,1,0}&$[3^{(3)}]$\\
B$_2$\hfill\FDynkin{1,2,0,0}&$[5,1^{(4)}]$
\end{tabular}\\
\\
\\
\\
\begin{tabular}[t]{rll}
\multicolumn3{l}{E$_6$, even subalgebra A$_5+$ A$_1$:}\\
\\
&nilpotent $e$&partitions in A$_5$, A$_1$\\
\hline
&A$_1$\hfill\EDynkin{1,0,0,0,0,0}{1}&$[1^{(6)}]$, $[2]$\bigstrut[t]\\
*&3A$_1$\hfill\EDynkin{0,0,0,1,0,0}{1}&$[2^{(3)}]$, $[1^{(2)}]$\\
&A$_2+$ A$_1$\hfill\EDynkin{1,1,0,0,0,1}{1}&$[3,1^{(3)}]$, $[2]$\\
*&2A$_2+$ A$_1$\hfill\EDynkin{0,1,0,1,0,1}{1}&$[3^{(2)}]$, $[2]$\\
&A$_3+$ A$_1$\hfill\EDynkin{1,0,1,0,1,0}{1}&$[4,2]$, $[1^{(2)}]$\\
&A$_4+$ A$_1$\hfill\EDynkin{1,1,1,0,1,1}{1}&$[5,1]$, $[2]$\\
&A$_5$\hfill\EDynkin{1,2,1,0,1,2}{1}&$[6]$, $[1^{(2)}]$
\end{tabular}&
\begin{tabular}[t]{ll}
\multicolumn2{l}{E$_6$, even subalgebra D$_5+T^1$:}\\
\\
nilpotent $e$&partition in D$_5$\\
\hline
2A$_1$\hfill\EDynkin{0,1,0,0,0,1}{1}&$[3,1^{(7)}]$\bigstrut[t]\\
A$_2+2$A$_1$\hfill\EDynkin{0,0,1,0,1,0}{1}&$[3^{(3)},1]$\\
A$_3$\hfill\EDynkin{2,1,0,0,0,1}{1}&$[5,1^{(5)}]$\\
D$_5(a_1)$\hfill\EDynkin{2,1,1,0,1,1}{1}&$[7,3]$
\end{tabular}
\end{tabular}%
}

\begin{tabular}{rll}
\multicolumn3{l}{E$_7$, even subalgebra D$_6+$ A$_1$:}\\
\\
&nilpotent $e$&partitions in D$_6$, A$_1$\\
\hline
&A$_1$\hfill\EDynkin{0,1,0,0,0,0,0}{1}&$[1^{(6)}]$, $[2]$\bigstrut[t]\\
&2A$_1$\hfill\EDynkin{0,0,0,0,0,1,0}{1}&$[3,1^{(9)}]$, $[1^{(2)}]$\\
*&$[3\rA_1]'$\hfill\EDynkin{0,0,1,0,0,0,0}{1}&$[2^{(6)}]$, $[1^{(2)}]$\\
*&4A$_1$\hfill\EDynkin{1,0,0,0,0,0,1}{1}&$[2^{(6)}]$, $[2]$\\
&A$_2+$ A$_1$\hfill\EDynkin{0,1,0,0,0,1,0}{1}&$[3^{(2)},1^{(6)}]$, $[2]$\\
&A$_2+2$A$_1$\hfill\EDynkin{0,0,0,1,0,0,0}{1}&$[3^{(3)},1^{(3)}]$, $[1^{(2)}]$\\
&A$_3$\hfill\EDynkin{0,2,0,0,0,1,0}{1}&$[7,5]$, $[1^{(2)}]$\\
*&2A$_2+$ A$_1$\hfill\EDynkin{0,0,1,0,0,1,0}{1}&$[3^{(4)}]$, $[2]$\\
&$[\rA_3+\rA_1]'$\hfill\EDynkin{0,1,0,1,0,0,0}{1}&$[4^{(2)},2^{(2)}]$, $[1^{(2)}]$\\
&A$_3+2$A$_1$\hfill\EDynkin{0,1,0,0,1,0,1}{1}&$[4^{(2)},2^{(2)}]$, $[2]$\\
&D$_4(a_1)+$ A$_1$\hfill\EDynkin{1,0,1,0,0,0,1}{1}&$[5,3,1^{(4)}]$, $[2]$\\
&A$_3+$ A$_2$\hfill\EDynkin{0,0,0,1,0,1,0}{1}&$[5,3^{(2)},1]$, $[1^{(2)}]$\\
&D$_4+$ A$_1$\hfill\EDynkin{1,2,1,0,0,0,1}{1}&$[7,5]$, $[2]$\\
&A$_4+$ A$_1$\hfill\EDynkin{0,1,0,1,0,1,0}{1}&$[5^{(2)},1^{(2)}]$, $[2]$\\
&D$_5(a_1)$\hfill\EDynkin{0,2,0,1,0,1,0}{1}&$[7,3,1^{(2)}]$, $[1^{(2)}]$\\
&$[\rA_5]'$\hfill\EDynkin{0,1,0,1,0,2,0}{1}&$[6^{(2)}]$, $[1^{(2)}]$\\
&$\rA_5+\rA_1$\hfill\EDynkin{0,1,0,1,0,1,2}{1}&$[6^{(2)}]$, $[2]$\\
&D$_6(a_2)$\hfill\EDynkin{1,0,1,0,1,0,2}{1}&$[7,5]$, $[1^{(2)}]$\\
&D$_5+$ A$_1$\hfill\EDynkin{1,2,1,0,1,1,0}{1}&$[9,3]$, $[2]$\\
&D$_6(a_1)$\hfill\EDynkin{1,2,1,0,1,0,2}{1}&$[9,3]$, $[1^{(2)}]$\\
&D$_6$\hfill\EDynkin{1,2,1,0,1,2,2}{1}&$[11,1]$, $[1^{(2)}]$
\end{tabular}

\resizebox{.9\textwidth}{!}{
\begin{tabular}[t]{rl|lr}
\multicolumn3{l}{E$_8$, even subalgebra E$_7+$ A$_1$:}\\
\\
&nilpotent $e$&label in E$_7$\hfill diagram in E$_7$& partition in A$_1$\\
\hline
&A$_1$\hfill\EDynkin{0,0,0,0,0,0,0,1}{1}&0&$[2]$\bigstrut[t]\\
*&3A$_1$\hfill\EDynkin{0,0,0,0,0,0,1,0}{1}&$[3\rA_1]''$\hfill\EDynkin{0,0,0,0,0,0,2}{1}&$[1^{(2)}]$\\
&A$_2+$ A$_1$\hfill\EDynkin{0,1,0,0,0,0,0,1}{1}&A$_2$\hfill\EDynkin{0,2,0,0,0,0,0}{1}&$[2]$\\
&A$_2+3$A$_1$\hfill\EDynkin{0,0,1,0,0,0,0,0}{1}&A$_2+3$A$_1$\hfill\EDynkin{2,0,0,0,0,0,0}{1}&$[1^{(2)}]$\\
*&2A$_2+$ A$_1$\hfill\EDynkin{0,1,0,0,0,0,1,0}{1}&2A$_2$\hfill\EDynkin{0,0,0,0,0,2,0}{1}&$[2]$\\
&A$_3+$ A$_1$\hfill\EDynkin{0,0,0,0,0,1,0,1}{1}&$[\rA_3+\rA_1]''$\hfill\EDynkin{0,2,0,0,0,0,2}{1}&$[1^{(2)}]$\\
&D$_4(a_1)+$ A$_1$\hfill\EDynkin{1,0,0,0,0,0,1,0}{1}&D$_4(a_1)$\hfill\EDynkin{0,0,2,0,0,0,0}{1}&$[2]$\\
&A$_3+$ A$_2+$ A$_1$\hfill\EDynkin{0,0,0,1,0,0,0,0}{1}&A$_3+$ A$_2+$ A$_1$\hfill\EDynkin{0,0,0,0,2,0,0}{1}&$[1^{(2)}]$\\
&D$_4+$ A$_1$\hfill\EDynkin{1,0,0,0,0,0,1,2}{1}&D$_4$\hfill\EDynkin{0,2,2,0,0,0,0}{1}&$[2]$\\
&A$_4+$ A$_1$\hfill\EDynkin{0,1,0,0,0,1,0,1}{1}&A$_4$\hfill\EDynkin{0,2,0,0,0,2,0}{1}&$[2]$\\
&A$_5$\hfill\EDynkin{0,2,0,0,0,1,0,1}{1}&$[\rA_5]''$\hfill\EDynkin{0,2,0,0,0,2,2}{1}&$[1^{(2)}]$\\
&D$_5(a_1)+$ A$_1$\hfill\EDynkin{0,0,0,1,0,0,0,2}{1}&D$_5(a_1)+$ A$_1$\hfill\EDynkin{0,2,0,0,2,0,0}{1}&$[1^{(2)}]$\\
&A$_4+$ A$_2+$ A$_1$\hfill\EDynkin{0,0,1,0,0,1,0,0}{1}&A$_4+$ A$_2$\hfill\EDynkin{0,0,0,2,0,0,0}{1}&$[2]$\\
&E$_6(a_3)+$ A$_1$\hfill\EDynkin{0,1,0,0,1,0,1,0}{1}&E$_6(a_3)$\hfill\EDynkin{0,0,2,0,0,2,0}{1}&$[2]$\\
&E$_7(a_5)$\hfill\EDynkin{0,0,0,1,0,1,0,0}{1}&E$_7(a_5)$\hfill\EDynkin{0,0,0,2,0,0,2}{1}&$[1^{(2)}]$\\
&D$_5+$ A$_1$\hfill\EDynkin{0,1,0,0,1,0,1,2}{1}&D$_5$\hfill\EDynkin{0,2,2,0,0,2,0}{1}&$[2]$\\
&A$_6+$ A$_1$\hfill\EDynkin{0,1,0,1,0,1,0,0}{1}&A$_6$\hfill\EDynkin{0,0,0,2,0,2,0}{1}&$[2]$\\
&E$_7(a_4)$\hfill\EDynkin{0,0,0,1,0,1,0,2}{1}&E$_7(a_4)$\hfill\EDynkin{0,2,0,2,0,0,2}{1}&$[1^{(2)}]$\\
&E$_6(a_1)+$ A$_1$\hfill\EDynkin{0,1,0,1,0,1,0,2}{1}&E$_6(a_1)$\hfill\EDynkin{0,2,0,2,0,2,0}{1}&$[2]$\\
&E$_7(a_3)$\hfill\EDynkin{0,2,0,1,0,1,0,2}{1}&E$_7(a_3)$\hfill\EDynkin{0,2,0,2,0,2,2}{1}&$[1^{(2)}]$\\
&E$_6+$ A$_1$\hfill\EDynkin{0,1,0,1,0,1,2,2}{1}&E$_6$\hfill\EDynkin{0,2,2,2,0,2,0}{1}&$[2]$\\
&E$_7(a_2)$\hfill\EDynkin{1,0,1,0,1,0,2,2}{1}&E$_7(a_2)$\hfill\EDynkin{2,2,2,0,2,0,2}{1}&$[1^{(2)}]$\\
&E$_7(a_1)$\hfill\EDynkin{1,2,1,0,1,0,2,2}{1}&E$_7(a_1)$\hfill\EDynkin{2,2,2,0,2,2,2}{1}&$[1^{(2)}]$\\
&E$_7$\hfill\EDynkin{1,2,1,0,1,2,2,2}{1}&E$_7$\hfill\EDynkin{2,2,2,2,2,2,2}{1}&$[1^{(2)}]$
\end{tabular}
}

\begin{tabular}[t]{rll}
\multicolumn3{l}{E$_8$, even subalgebra D$_8$:}\\
\\
&nilpotent $e$&partition in D$_8$\\
\hline
&2A$_1$\hfill\EDynkin{0,1,0,0,0,0,0,0}{1}&$[3,1^{(13)}]$\bigstrut[t]\\
*&4A$_1$\hfill\EDynkin{1,0,0,0,0,0,0,0}{1}&$[2^{(8)}]$\\
&A$_2+2$A$_1$\hfill\EDynkin{0,0,0,0,0,1,0,0}{1}&$[3^{(3)},1^{(7)}]$\\
&A$_3$\hfill\EDynkin{0,1,0,0,0,0,0,2}{1}&$[5,1^{(11)}]$\\
*&2A$_2+2$A$_1$\hfill\EDynkin{0,0,0,0,1,0,0,0}{1}&$[3^{(5)},1]$\\
&A$_3+2$A$_1$\hfill\EDynkin{0,0,1,0,0,0,0,1}{1}&$[4^{(2)},2^{(4)}]$\\
&A$_3+$ A$_2$\hfill\EDynkin{0,1,0,0,0,1,0,0}{1}&$[5,3^{(2)},1^{(5)}]$\\
*&2A$_3$\hfill\EDynkin{0,1,0,0,1,0,0,0}{1}&$[4^{(4)}]$\\
&D$_5(a_1)$\hfill\EDynkin{0,1,0,0,0,1,0,2}{1}&$[7,3,1^{(6)}]$\\
&A$_4+2$A$_1$\hfill\EDynkin{0,0,0,1,0,0,0,1}{1}&$[5^{(2)},3,1^{(3)}]$\\
*&A$_4+$ A$_3$\hfill\EDynkin{0,0,0,1,0,0,1,0}{1}&$[5^{(3)},1]$\\
&A$_5+$ A$_1$\hfill\EDynkin{0,1,0,1,0,0,0,1}{1}&$[6^{(2)},2^{(2)}]$\\
&D$_5(a_1)+$ A$_2$\hfill\EDynkin{0,0,1,0,0,1,0,1}{1}&$[7,3^{(3)}]$\\
&D$_6(a_2)$\hfill\EDynkin{1,0,1,0,0,0,1,0}{1}&$[7,5,1^{(4)}]$\\
&D$_6(a_1)$\hfill\EDynkin{1,0,1,0,0,0,1,2}{1}&$[9,3,1^{(4)}]$\\
&D$_6$\hfill\EDynkin{1,2,1,0,0,0,1,2}{1}&$[11,1^{(5)}]$\\
&D$_7(a_2)$\hfill\EDynkin{0,1,0,1,0,1,0,1}{1}&$[9,5,1^{(2)}]$\\
*&A$_7$\hfill\EDynkin{0,1,0,1,0,1,1,0}{1}&$[8^{(2)}]$\\
&D$_7$\hfill\EDynkin{1,2,1,0,1,1,0,1}{1}&$[13,1^{(3)}]$
\end{tabular}

\

\

\subsubsection{\bf Remark}
Note that not all possible fixed point algebras of involutive automorphisms are realized as $\g^{\mathrm{ev}}$ for some nilpotent element. Indeed, the subalgebra $\g^{\mathrm{ev}}$ of $\g$ is the fixed point set of an involutive automorphism of $\g$, which lies in the center of the subgroup SL(2) of $G$
with Lie algebra $\sa$, acts as 1 on $\g^{\mathrm{ev}}$ and as $-1$ on the odd part of the grading.
This rules out some of the fixed point subalgebras, listed above, as $\g^{\mathrm{ev}}$. For example, this rules out $\g^{\mathrm{ev}}$ in E$_7$ of types E$_6+T^1$ and
A$_7$. All other possibilities in exceptional Lie algebras do occur. For classical types, all possibilities are realized for type A, all semisimple $\g^{\mathrm{ev}}$ occur for types B, C, D, and, in addition, the subalgebra D$_{2m}+T^1$ occurs for D$_{2m+1}$.

\

\

\

\subsection{Algebra $(\g_{-d},*)$ for mixed type nilpotent elements}\label{appB}

Here we describe the algebra structures $(\g_{-d},*)$ for nilpotent elements $e$ of mixed type.

Let us recall from \cite{ekv}*{Remark 3.2} that reducing subalgebras $\qa$ for such $e$ can be defined as semisimple subalgebras normalized by the $\sla(2)$-triple $\sa$ for $e$ such that in the decomposition $e=e_\qa+e'$, where $e_\qa\in\qa$, $e'\in\z(\qa)$, the nilpotent element $e_\qa$ has the same depth and rank in $\qa$ as $e$ in $\g$. We then have

\subsubsection{\bf Proposition}
\emph{Let $\qa$ be a reducing subalgebra in the above sense, for any $e\in\g$ (of even depth). Then for any $F,F'\in\qa_{-d}$ their $*$-product in $\g_{-d}$ induced by $e$ coincides with that induced by $e_\qa$. In particular, $\qa_{-d}\subseteq\g_{-d}$ is a $*$-subalgebra.}
\begin{proof}
From $e=e_\qa+e'$ with $e'\in\z(\qa)$, it follows that $(\ad e)x=(\ad e_\qa)x$ for any $x\in\qa$. Thus for $F,F'\in\qa_{-d}$ we have
$$
[(\ad e)^{\frac d2}F,F']=[(\ad e_\qa)^{\frac d2}F,F'],
$$
i. e. the two $*$-products on $\qa_{-d}$ coincide.
\end{proof}

Moreover it is shown in \cite{ekv} that for any $e$ of mixed type there is a reducing subalgebra $\qa$ in this sense such that $e_\qa$ is of semisimple type in $\qa$.

This is used in \cite{ekv} to group nilpotent elements into \emph{bushes}; each bush is a subset of nilpotent elements admitting a common reducing subalgebra $\qa$ with the same $e_\qa$, the latter being the unique nilpotent element of semisimple type in the bush.

In particular, if $\qa_{-d}=\g_{-d}$ then the *-algebra structure on $\g_{-d}$ is one of those corresponding to a nilpotent element of semisimple type that we have already described. It thus remains to consider the cases when for any reducing subalgebra $\qa$ with $e_\qa$ of semisimple type in $\qa$, the space $\qa_{-d}$ is a proper subalgebra of $\g_{-d}$.

Note that such $e$ can be also characterized using the particular reducing subalgebra described in \cite{ekv}*{Proposition 3.10}: these are precisely the nilpotent elements with the property that, for the $\sla(2)$-triple $\sa$ of $e$ in the reducing subalgebra $\qa$ generated by the $\sa$-submodule of $\g$ generated by $\g_{-d}$, $e_\qa$ is not of semisimple type in $\qa$.

In what follows we will encounter commutative algebras over $\field$ of the following kind.

We will denote by $J_n^c(BD)$, $c\in\field$, the commutative algebra of dimension $n+1$, with basis $\bf1$, $x_1$, ..., $x_n$ and multiplication table
$$
{\bf1}x_i=x_i,\quad x_i^2={\bf1},\quad x_ix_j=0 \text{ for $i\ne j$, and } {\bf1}^2=c{\bf1}.
$$

Furthermore, let $(H_8,*)$ denote the 8-dimensional space of traceless $3\times 3$ matrices, with the multiplication
$$
A*B:=\frac{AB+BA}2-\frac13\operatorname{Trace}\left(\frac{AB+BA}2\right),
$$
and let $H_5\subset H_8$ be its 5-dimensional subspace consisting of symmetric matrices. Clearly then $H_5$ is a $*$-subalgebra of $H_8$. It contains the subalgebra of diagonal matrices isomorphic to $\CA_{-1}(2)$, as well as infinitely many subalgebras isomorphic to $J^{-1}_2(BD)$, for example the subalgebra spanned by diagonal matrices and any one of the $e_{12}+e_{21}$, $e_{13}+e_{31}$, or $e_{23}+e_{32}$ is such.

Thus $J_n^1(BD)$ is isomorphic to the Jordan algebra $J_n(BD)$. Moreover a calculation, similar to that in Lemma \ref{case4}, shows that the algebra $J_1^c(BD)$ is isomorphic to $\CA_{\frac{c-1}2}(2)$. For most other values of $n$ and $c$ this algebra does not have unity and is not Jordan, neither does it satisfy the identities \eqref{confass} or \eqref{mockjordan}. Note that $J_n^c(BD)$ contains isomorphic copies of $J_m^c(BD)$ for $m\le n$.

Note also that the $*$-multiplication on $H_8$ is the unique commutative multiplication invariant under the adjoint action of $\sla(3)$ on it, while the $*$-multiplication on $H_5$ is the unique commutative multiplication invariant under the action of $\so(3)\cong\sla(2)$ realizing $H_5$ as the $5$-dimensional irreducible representation of $\sla(2)$ ($=$ the $5$-dimensional irreducible summand of the symmetric square of the adjoint representation of $\so(3)$).

For classical type Lie algebras $\g$, we have the following cases when $\g_{-d}$ is strictly larger than the $-d$ degree component for the nilpotent element of the semisimple type in the same bush:

In $\so((2k+1)(2\ell+1)+n_1+...+n_j)$, the nilpotent element $e$ with the orbit partition $[(2k+1)^{(2\ell+1)},n_1,...,n_j]$, $k,\ell\ge1$, $n_1,...,n_j<2k+1$ --- depth is $4k$, with $e_\qa$ having partition $[(2k+1)^{2\ell}]$ in the reducing subalgebra $\qa=\so(2\ell(2k+1))$. Then the algebra $(\g_{-d},*)$ is $\so(2\ell+1)$ with the adjoint action of $Z(\sa)$, while in the reducing subalgebra $\qa$ its subalgebra $\qa_{-d}$ is isomorphic to $\so(2\ell)$.

In $\so((2k+1)(\ell+1)+2+n_1+...+n_j)$, the nilpotent element $e$ with the orbit partition $[2k+3,(2k+1)^{(\ell)},n_1,...,n_j]$, $k,\ell\ge1$, $n_1,...,n_j<2k+1$ --- depth is $4k+2$, with $e_\qa$ having partition $[2k+3,2k+1]$ in the reducing subalgebra $\qa=\so(4k+4)$. Here the algebra $(\g_{-d},*)$ is isomorphic to $J^{-(2k-1)}_\ell(BD)$. Its subalgebra $\qa_{-d}\subseteq\g_{-d}$ is isomorphic to $J^{-(2k-1)}_1(BD)\cong\CA_{-k}(2)$.

It follows from the description of bushes for algebras of classical types in \cite{ekv}*{end of Section 4} that the above are the only cases for classical types when $\dim\g_{-d}$ is larger than that for the element of the semisimple type in the bush.

For exceptional type Lie algebras $\g$, nilpotent elements $e$ such that for any reducing subalgebra $\qa$ with $e_\qa$ of semisimple type one has $\qa_{-d}\subsetneqq\g_{-d}$ are the following:

F$_4$, label ${\mathrm A}_2+\tilde{\mathrm A}_1$: depth is 4, the algebra $(\g_{-d},*)$ is isomorphic to $\sla(2)$, realizing the adjoint representation of $\z(\sa)\cong\sla(2)$. The subalgebra $\qa$ of $\g$ generated by the $\sa$-submodule of $\g$ generated by the 1-dimensional Cartan subalgebra $(\qa_{-d},0)$ of $(\g_{-d},*)$ is of type A$_2$, and in the decomposition $e=e_\qa+e'$ of $e$ in $\qa\oplus\z(\qa)$ the nilpotent element $e_\qa$ is principal in $\qa$. It has label A$_2$ in $\g$ and $(e,e_\qa)$ constitute a bush in F$_4$.

E$_6$, label ${\mathrm A}_2+2{\mathrm A}_1$: depth is 4, the algebra $(\g_{-d},*)$ is isomorphic to $\sla(2)$, realizing the adjoint representation of $\z(\sa)\cong\sla(2)$. The subalgebra $\qa$ of $\g$ generated by the $\sa$-submodule of $\g$ generated by the Cartan subalgebra of $(\g_{-d},*)$ is of type A$_2$, and in the decomposition $e=e_\qa+e'$ of $e$ in $\qa\oplus\z(\qa)$ the nilpotent element $e_\qa$ is principal in $\qa$. It has label A$_2$ in $\g$ and $(e,e_\qa)$ together with the nilpotent element with label ${\mathrm A}_2+{\mathrm A}_1$ (having $\dim(\g_{-d})=1$) constitute a bush in E$_6$.

E$_7$, label ${\mathrm A}_2+2{\mathrm A}_1$: same properties as the element with the same label in E$_6$, except that the bush contains one more element, with label ${\mathrm A}_2+3{\mathrm A}_1$ (see next entry).

E$_7$, label ${\mathrm A}_2+3{\mathrm A}_1$: depth 4, the algebra $(\g_{-d},*)$ is isomorphic to the simple Malcev algebra $M$ of dimension $7$, realizing the smallest irreducible representation of $\z(\sa)$, which is of type G$_2$. The subalgebra $\qa$ of $\g$ generated by the $\sa$-submodule of $\g$ generated by the Cartan subalgebra of $(\g_{-d},*)$ is of type A$_2$, with $e_\qa$ principal there. Moreover $(\g_{-d},*)$ admits an infinite family of 3-dimensional subalgebras, each isomorphic to $\sla(2)$. For the reducing subalgebras $\qa$ generated by the $\sa$-submodules generated by any one of those $\sla(2)$-subalgebras of $(\g_{-d},*)$, the element $e_\qa$ has label ${\mathrm A}_2+2{\mathrm A}_1$ in $\g$.

E$_7$, label ${\mathrm A}_4+{\mathrm A}_2$: depth is 8, $(\g_{-d},*)$ is isomorphic to $\sla(2)$, realizing the adjoint representation of $\z(\sa)\cong\sla(2)$. For the reducing subalgebra $\qa$ generated by the $\sa$-submodule generated by the Cartan subalgebra of $(\g_{-d},*)$, the nilpotent element $e_\qa$ is of semisimple type; in $\g$ it has label A$_4$. The bush also contains the nilpotent element with label ${\mathrm A}_4+{\mathrm A}_1$, with $\dim\g_{-d}=1$.

E$_7$, label ${\mathrm A}_3+{\mathrm A}_2$: depth 6, the $3$-dimensional algebra $(\g_{-d},*)$ is isomorphic to $J^{-1}_2(BD)$. Here $\z(\sa)$ is a 1-dimensional torus acting on $\g_{-d}$ with eigenvalues $\pm1$ and $0$. To obtain the element of semisimple type from the bush we may take any subalgebra of $J^{-1}_2(BD)$ spanned by $\bf1$ and some element $x$ with $x*x=\bf1$. This subalgebra is isomorphic to $\CA_{-1}(2)$ and the $\sa$-submodule generated by it generates a reducing subalgebra $\qa$ such that $e_\qa$ is of semisimple type in it. In $\g$ it has label D$_4(a_1)$. The bush also contains an element with label D$_4(a_1)+$A$_1$, with $\g_{-d}$ the same as for $e_\qa$, as well as one more element (see the next entry).

E$_7$, label ${\mathrm A}_3+{\mathrm A}_2+{\mathrm A}_1$: depth 6, the algebra $(\g_{-d},*)$ is isomorphic to $H_5$. Here $\z(\sa)$ is $\sla(2)$, and its representation on $\g_{-d}$ is irreducible. 3-dimensional subalgebras of $H_5$ isomorphic to $J^{-1}_2(BD)$ realize, by the same procedure, nilpotent elements with label ${\mathrm A}_3+{\mathrm A}_2$.

E$_8$, labels ${\mathrm A}_2+2{\mathrm A}_1$ and ${\mathrm A}_2+3{\mathrm A}_1$ --- this bush has exactly the same properties as the one with these labels in E$_7$.

E$_8$, label ${\mathrm A}_4+{\mathrm A}_2$ --- same as the nilpotent element with this label in E$_7$, but the bush contains two more elements: the one with label ${\mathrm A}_4+2{\mathrm A}_1$, with $\dim\g_{-d}=1$, and the one described in the next entry.

E$_8$, label ${\mathrm A}_4+{\mathrm A}_2+{\mathrm A}_1$ --- depth is 8 and $\dim\g_{-d}=3$; the algebra $(\g_{-d},*)$ is the same as the one for the element with label ${\mathrm A}_4+{\mathrm A}_2$ in the same bush.

E$_8$, labels ${\mathrm A}_3+{\mathrm A}_2$ and ${\mathrm A}_3+{\mathrm A}_2+{\mathrm A}_1$ --- same properties as the ones of this bush in E$_7$, but the bush here contains one more element, see the last entry.

E$_8$, label D$_4(a_1)+{\mathrm A}_2$: here, as for other elements in the bush, depth is 6. The algebra $(\g_{-d},*)$ is isomorphic to $H_8$, realizing the adjoint action of $\z(\sa)$ which in this case is $\sla(3)$. The algebra $H_8$ contains infinitely many 5-dimensional subalgebras giving rise to nilpotent elements with label ${\mathrm A}_3+{\mathrm A}_2+{\mathrm A}_1$ from the bush. For example, $H_5$ is such, but also isomorphic to $H_5$ is the subalgebra of $H_8$ spanned by the diagonals, two of the antisymmetric matrices $e_{12}-e_{21}$, $e_{13}-e_{31}$, $e_{23}-e_{32}$ and their $*$-product, which is symmetric, e.~g. $e_{12}-e_{21}$, $e_{23}-e_{32}$ and $e_{13}+e_{31}$.

\subsection{Chains of nilpotent elements}

Recall \cite{ekv} that any nilpotent element $e\in\g$ not of nilpotent type uniquely decomposes in a sum of commuting elements: $e=e^s+e^n$, where $e^s$ lies in the minimal reducing subalgebra $\qa_{\min}$ and $e^n$ lies in its centralizer. The nilpotent element $e^s$ is of semisimple type in $\g$, and $e^n$ can be of any type. Let $f(e)=e^n$. Thus $f(e)=0$ for $e$ of semisimple type; for $e$ of nilpotent type it is natural to put $f(e)=e$. Then for each nilpotent element $e$ we have a \emph{chain}
$$
e=f^0(e),f(e)=f^1(e),f(f(e))=f^2(e),...,f^{\ell-1}(e),
$$
where $\ell$, the \emph{length} of the chain for $e$, is the smallest natural number such that the iterate $f^{\ell-1}(e)$ is of either semisimple or of nilpotent type. Thus for $e$ of mixed type $\ell\ge2$.

If $\g$ is of classical type, and $e$ is a nilpotent element, corresponding to the partition $(n_1,...,n_k)$, with $n_1\ge\cdots\ge n_k>0$, then $e^s$ corresponds to the partition $(n_1,...,n_j,1^{(n_{j+1}+...+n_k)})$ for some $0\le j\le k$, and $e^n=f(e)$ corresponds to the partition $(n_{j+1},...,n_k,1^{(n_1+...+n_j)})$. According to  \cite{ekv}*{p. 111}, except for $e$ of nilpotent type, here $j$ is the largest natural number with the property that $(n_1,...,n_j,1^{(n_{j+1}+...+n_k)})$ is the partition of a nilpotent of semisimple type in $\g$. This rule determines the chain for $e$.

For orthogonal Lie algebras the chain can terminate with an element of nilpotent type. One can show that this happens if and only if in the corresponding partition $n_1\ge\cdots\ge n_k$, there is an odd $n_j$ with $n_{j+1}=n_j-1$ such that the maximal subsequence $n_i\geqslant\cdots\geqslant n_j$ consisting of consecutive odd numbers (repetitions allowed) has odd sum.

\subsubsection{\bf Examples}
In B$_{27}$, there is a chain
$$
(7,5^{(4)},4^{(4)},2^{(2)},1^{(8)})\mapsto(5^{(3)},4^{(4)},2^{(2)},1^{(20)})\mapsto(5,4^{(4)},2^{(2)},1^{(30)}),
$$
the last one is of nilpotent type.

In D$_9$,
$$
(5^{(3)},3)\mapsto(5,3,1^{(10)}),
$$
the last one is of semisimple type.

In C$_{17}$
$$
(5^{(2)},4^{(3)},3^{(4)})\mapsto(4^{(3)},3^{(4)},1^{(10)})\mapsto(3^{(4)},1^{(22)}),
$$
the last one is of semisimple type.

\

For $\g$ of exceptional types, the length of all mixed type nilpotent elements $e$ is equal to $2$, with two exceptions, both in $\rE_8$, when the length is $3$:
$$
\rA_4+\rA_2+\rA_1\mapsto\rA_2+\rA_1\mapsto\rA_1\text{ and }\rD_5(a_1)+\rA_2\mapsto\rA_2+2\rA_1\mapsto2\rA_1.
$$
Moreover, for exceptional types all ending elements of chains for mixed types are of semisimple type, with one exception, again in $\rE_8$,
which is the last entry for $\rE_8$ below.

The chains of length 2 for mixed type nilpotent elements in $\g$ of exceptional type are as follows:

\

\noindent In $\rE_6$,

all of $\rA_2+\rA_1$, $\rA_3+\rA_1$, $\rA_4+\rA_1$ go to $\rA_1$ in one step;

$\rA_2+2\rA_1$, $\rD_5(a_1)$ go to $2\rA_1$.

\noindent In $\rE_7$,

$\rA_2+\rA_1$, $[\rA_3+\rA_1]'$, $[\rA_3+\rA_1]''$, $\rA_4+\rA_1$, $\rA_5+\rA_1$, $\rD_4+\rA_1$, $\rD_4(a_1)+\rA_1$, $\rD_5+\rA_1$,

and $\rE_7(a_3)$ go to $\rA_1$;

$\rA_2+2\rA_1$, $\rA_3+2\rA_1$, $\rA_3+\rA_2$, $\rD_5(a_1)$, $\rD_6(a_1)$ go to $2\rA_1$;

$\rA_2+3\rA_1$, $\rA_3+\rA_2+\rA_1$, $\rD_5(a_1)+\rA_1$, $\rE_7(a_2)$, $\rE_7(a_4)$ go to $[3\rA_1]''$;

$\rA_4+\rA_2$ goes to $\rA_2$.

\noindent In $\rE_8$,

$\rA_2+\rA_1$, $\rA_3+\rA_1$, $\rD_4(a_1)+\rA_1$, $\rD_4+\rA_1$, $\rA_4+\rA_1$, $\rA_5+\rA_1$, $\rE_6(a_3)+\rA_1$, $\rD_5+\rA_1$,

$\rA_6+\rA_1$, $\rE_6(a_1)+\rA_1$, $\rE_7(a_3)$, $\rE_6+\rA_1$, $\rE_8(b_4)$, $\rE_8(a_3)$ go to $\rA_1$;

$\rD_4(a_1)+\rA_2$, $\rA_4+\rA_2$, $\rD_4+\rA_2$, $\rD_5+\rA_2$, $\rE_8(b_6)$, $\rE_8(b_5)$ go to $\rA_2$;

$\rD_7\left(a_2\right)$ goes to $\rA_3$;

$\rA_2+2\rA_1$, $\rA_3+2\rA_1$, $\rA_3+\rA_2$, $\rD_5(a_1)$, $\rA_4+2\rA_1$, $\rD_6(a_1)$, $\rD_7(a_1)$ go to $2\rA_1$;

$\rA_2+3\rA_1$, $\rA_3+\rA_2+\rA_1$, $\rD_5(a_1)+\rA_1$, $\rE_7(a_4)$, $\rE_7(a_2)$ go to $3\rA_1$.

\noindent In $\rF_4$,

$\rC_3(a_1)\mapsto\rA_1$,

$\rA_2+\tilde\rA_1\mapsto\tilde\rA_1$.

\end{document}